\documentclass[10pt]{elsarticle}
\usepackage[margin=1.5in]{geometry}
\usepackage{amssymb, amsmath, mathrsfs}
\usepackage{bm}
\usepackage{latexsym}
\usepackage{graphicx}
\usepackage{subfig, float}
\usepackage{caption}
\usepackage[colorlinks]{hyperref}
\usepackage{tabularx}
\usepackage{multirow}
\usepackage{lineno}
\graphicspath{{./fig/}}
\usepackage{color, xcolor}

\def\D{\mathbb{D}}
\def\I{\mathbb{I}}
\def\S{\mathbb{S}}

\newtheorem{remark}{Remark}

\newtheorem{theorem}{Theorem}
\allowdisplaybreaks[4]
\begin{document}

\begin{frontmatter}
    \title{Modeling and simulation of inductionless magnetohydrodynamic free surface problems with unmatched densities}
    \author{Jiancheng Wang}
    \ead{202311110403@std.uestc.edu.cn}
    \author{Maojun Li}
    \ead{limj@uestc.edu.cn}
    \author{Zeyu Xia}
    \ead{zeyuxia@uestc.edu.cn}
    \author{Liwei Xu}
    \ead{xul@uestc.edu.cn}
    \address{School of Mathematical Sciences, University of Electronic Science and Technology of China, Sichuan, 611731, P.R. China }

    \journal{***}

    \begin{abstract}
        We propose a new diffuse interface model for simulating an inductionless magnetohydrodynamic (MHD) free surface problem. By using the Onsager's variational principle and the laws of thermodynamics, we derive a thermodynamically consistent system that couples the Cahn--Hilliard equation modeling phase separation, the Navier--Stokes equations governing fluid motion, and a generalized Darcy's law accounting for electromagnetic effects. In contrast to existing diffuse interface MHD models, the proposed model can handle general material properties in practical engineering applications. Furthermore, through asymptotic arguments, we investigate the sharp interface limit, and then demonstrate that the classical sharp interface model can be recovered as the interface thickness approaches zero, theoretically validating the proposed diffuse interface model as an approximate approach. An efficient decoupled, linear, and charge-conservative finite element scheme is designed, and it significantly facilitates the large-scale and accurate numerical simulations involving large parameter ratios. Finally, we present several three-dimensional numerical experiments of magnetic damping effects on bubble dynamics for the demonstration of the capability of the proposed model and method in capturing complex MHD phenomena.
    \end{abstract}


\end{frontmatter}


\section{Introduction}

When the magnetic fields are applied to the electrically conducting and nonmagnetic fluids (such as the liquid metals and strong electrolytes), the induced magnetohydrodynamic (MHD) effects will significantly influence the fluid dynamics, resulting in a fundamentally distinct flow behavior. Considering the ubiquitous presence of such a scenario in industrial processes, numerous studies have been carried out to explore the underlying flow mechanism, such as vortex dynamics, heat transfer, and turbulence, see, for example, \cite{2017_Davidson, 2025_Fan, 2018_Hiptmair}. While these flows exhibit extremely complex dynamics in practical engineering applications, free surface problems are usually  encountered, yielding a more intricate fluid motion. Representative examples include continuous casting and refining of metals \cite{2024_Gou, 2018_Thomas}, magnetic destabilization of free surfaces in aluminum reduction cells \cite{2003_Gerbeau, 2019_Herreman}, liquid metal batteries \cite{2018_Tucs}, and plasma-facing components \cite{2004_Morley}, among others.

In this work, we consider the inductionless MHD free surface problems in a fixed bounded domain $\Omega\subset\mathbb{R}^3$, which contains two incompressible and immiscible fluids occupying two open, time-dependent subdomains $\Omega_\pm(t)$, respectively. Characterized by the fact that the free surfaces are codimension-1 manifolds from a macroscopic perspective \cite{2011_Gross}, the two subdomains $\Omega_\pm(t)$ are separated by a sharp interface $\Gamma(t)=\overline{\Omega_+(t)}\cap\overline{\Omega_-(t)}$. Consequently, the flow dynamic behavior can be governed by the following classical sharp interface model:
\begin{subequations}
\label{sharp_interface}
    \begin{align}
         \rho_\pm\left(\partial_t\bm{u}+\bm{u}\cdot\nabla\bm{u}\right) & =\nabla\cdot(2\eta_{\pm}\D(\bm{u})-p\I)+\bm{J}\times\bm{B}+\rho_\pm\bm{g} & & \text{in }\Omega_\pm(t), \label{sharp_1} \\
        \nabla\cdot\bm{u} & =0 & & \text{in }\Omega_\pm(t), \label{sharp_2} \\
        \bm{J} & =\sigma_\pm\left(\bm{E}+\bm{u}\times\bm{B}\right) & & \text{in }\Omega_\pm(t), \label{sharp_3} \\
        \bm{E} & =-\nabla{V} & & \text{in }\Omega_\pm(t), \label{sharp_4} \\
        \nabla\cdot\bm{J} & =0 & & \text{in }\Omega_\pm(t), \label{sharp_5} \\
        [\![2\eta\D(\bm{u})-p\I]\!]\bm{n}_\Gamma & =-\lambda\kappa\bm{n}_\Gamma & & \text{on }\Gamma(t), \label{interface_condition_1} \\
        [\![\bm{u}]\!] & = {\bm 0} & & \text{on }\Gamma(t), \label{interface_condition_2} \\
        [\![\bm{J}\cdot\bm{n}_\Gamma]\!] & =0 & & \text{on }\Gamma(t), \label{interface_condition_33} \\
        [\![V]\!] & =0 & & \text{on }\Gamma(t), \label{interface_condition_4} \\
        {V}_\Gamma & =\bm{u}\cdot\bm{n}_\Gamma & & \text{on }\Gamma(t). \label{interface_condition_5}
    \end{align}
\end{subequations}
In this system, the unknowns are the fluid velocity $\bm{u}$, the pressure $p$, the current density $\bm{J}$, the electric field $\bm{E}$ and the electrostatic potential $V$. In the above equations, $\D(\bm{u})=\frac{1}{2}\big(\nabla\bm{u}+(\nabla\bm{u})^\top\big)$ and $\I$ stand for the deformation and the identity tensor, respectively. $\bm{g}=(0,0,-g)^{\top}$ represents the gravitational acceleration with $g$ being a positive constant. $\bm{B}$ is the uniform external magnetic field with the magnitude being denoted by $B_r$. The material properties $\rho_\pm$, $\eta_\pm$ and $\sigma_\pm$ are the constant densities, the constant dynamic viscosities, and the constant electrical conductivities of the two fluids, respectively. Moreover, we denote by $\lambda$ the surface tension, $\kappa$ the mean curvature, $\bm{n}_\Gamma$ the unit normal on $\Gamma(t)$ pointing to $\Omega_+(t)$, and $V_\Gamma$ the normal velocity of $\Gamma(t)$. We use the notation $\eta=\eta_+\mathscr{X}_{\Omega_+(t)}+\eta_-\mathscr{X}_{\Omega_-(t)}$ with $\mathscr{X}_{\Omega_\pm(t)}$ being the characteristic functions corresponding to $\Omega_\pm(t)$. Generally, $[\![\chi]\!]=\chi_+-\chi_-$ denotes the jump of quantity $\chi$ across $\Gamma(t)$ from $\Omega_+(t)$ to $\Omega_-(t)$. Equations \eqref{sharp_1}-\eqref{sharp_5} stand for the inductionless MHD equations, which can be employed to simulate the mutual interaction between the electromagnetic fields and the fluid flows in most laboratory experiments and industrial processes; see \cite{2017_Davidson} for a physical justification. Meanwhile, equations \eqref{interface_condition_1}-\eqref{interface_condition_4} are the jump conditions across the interface $\Gamma(t)$. Equation \eqref{interface_condition_5} models the movement of $\Gamma(t)$, accounting for the immiscibility assumption.

Numerically solving the system \eqref{sharp_interface} typically presents a significant challenge primarily due to the presence of the free surfaces and the jump conditions. To characterize the evolving interfaces of MHD flows, substantial numerical methodologies have been developed in the past decades. The level-set method \cite{2010_Flueck, 2004_Morley, 2006_Munger} was one of the earliest and extensively utilized methods for MHD free surface simulations. In this approach, the interface is represented by a signed distance function, which attains a zero value on the interface. Several numerical techniques are employed to achieve a stable simulation, such as the re-initialization (or re-distancing) method and the continuum surface force model; see \cite{2011_Gross} for more details of these techniques. Particularly, in \cite{2010_Flueck} the level-set algorithm only serves as an indicator of the interface, and an interpolation method was employed to construct a fitted mesh. Moreover, Zhang and Ni \cite{2014_Zhang, 2014_Zhang_pof, 2016_Zhang} employed the volume-of-fluid method to characterize the free surface. In this method, a volume fraction function is exploited to describe the interface, and some interface reconstruction techniques are required to maintain interfacial sharpness. Furthermore, moving mesh methods have been extensively applied in simulating MHD free surface problems. Gerbeau et al. \cite{2003_Gerbeau} utilized an arbitrary Lagrangian–Eulerian approach for the interface motion and developed a structure-preserving scheme. Samulyak et al. \cite{2007_Samulyak} employed a front-tracking technique for capturing the interface. Pan et al. \cite{2018_Pan} exploited the immersed boundary method to accurately resolve the moving boundary. Despite the universal and successful application of the above algorithms in MHD free surface simulations, extra numerical efforts are usually required to obtain the stable simulation results, including the high-resolution interface \cite{2011_Gross}.

In addition to the aforementioned approaches, the diffuse interface method \cite{1998_Anderson} emerges as another effective tool for free surface modeling. This methodology regards the interface as a significantly thin transition layer where the fluid properties vary smoothly. The derived models inherently incorporate the thermodynamic structure and the energy stability property. Meanwhile, this approach naturally handles the topological changes of the interface. Correspondingly, it requires solving a fourth-order nonlinear partial differential equation, i.e., the Cahn--Hilliard equation. This numerical challenge can be efficiently addressed by the mixed finite element methods \cite{2024_Wang} and some linearizing techniques \cite{2015_Shen}. Recently, Yang et al. \cite{2019_YangJ} proposed a diffuse interface model to simulate MHD problems and established the existence of the weak solutions. Subsequently, many studies have been devoted to the convergence analysis of numerical schemes for this model; see, e.g., \cite{2023_Qiu, 2025_Su, 2024_WangC, 2025_Yang}. Additionally, Chen and Zhang \cite{2020_Chen} proposed another diffuse interface MHD model, with its existence of the weak solutions being established by Wang et al. \cite{2025_Wang}. These models possess two features: (i) they are not derived from the thermodynamic principles rather than a direct coupling of the Cahn--Hilliard equation with different MHD equations; and (ii) they are not serving as the approximations of the corresponding sharp interface MHD models such as the system \eqref{sharp_interface}. Therefore, these features limit their application in the MHD simulations for general industrial processes.

This paper firstly aims to propose a novel diffuse interface model which are equivalent to the system \eqref{sharp_interface} under proper assumptions. The proposed model is rigorously derived from the thermodynamic principles \cite{2012_Abels} via the Onsager's variational principle \cite{2006_Qian}. As a result, the incompressibility constraint remains and an extra flux is incorporated into the momentum equation to guarantee the thermodynamic consistency, in turn leading to an energy law that is not an obvious approximation of the sharp interface counterpart, see Remark \ref{energy_relation} in Section \ref{section_2}. Utilizing an effective parameter scaling law proposed by Magaletti et al. \cite{2013_Magaletti}, and the method of formally matched asymptotic expansions \cite{2012_Abels, 2018_Xu}, we rigorously demonstrate that the system \eqref{sharp_interface} can be recovered as the interface thickness tends to zero. This result validates the reliability of the diffuse interface model in practical MHD simulations on a theoretical level. Secondly, we develop an efficient backward differential formula of order two (BDF2) finite element scheme, incorporating the linearizing technique \cite{2015_Shen} and the pressure stabilization method \cite{2009_Guermond} to achieve a decoupled and linear structure, to perform the numerical simulations. Meanwhile, such a scheme is charge-conservative (the current density is exactly solenoidal on the discrete level), in turn enabling an accurate MHD simulation \cite{2014_Zhang}. Focusing on applications in metallurgy processes, we compute several three-dimensional numerical examples to explore the magnetic damping effects on bubble dynamics, demonstrating the robustness and applicability of the proposed model and method in realistic MHD simulations.

The remainder of this paper is organized as follows. In Section \ref{section_2} we present the diffuse interface model and discuss the sharp interface asymptotics. Then, in Section \ref{section_3}, the discrete scheme is presented, and several numerical results are presented. Conclusions and future research directions are discussed in Section \ref{section_4}.

\section{Diffuse interface model}
\label{section_2}

\subsection{Governing equations}
Considering the approximation of the sharp interface $\Gamma(t)$ in the system \eqref{sharp_interface} as an artificial interface with thickness $\epsilon$, we propose the following diffuse interface model:
\begin{subequations}
\label{phase_field}
    \begin{align}
        \partial_t\phi+\nabla\cdot(\phi\bm{u}) & =\nabla\cdot\big(m_\epsilon(\phi)\nabla\mu\big), \label{phase_model_1} \\
        \mu & =-\widehat{\lambda}\epsilon\Delta\phi+\widehat{\lambda}\epsilon^{-1}(\phi^3-\phi), \label{phase_model_2} \\
        \rho(\phi)\partial_t\bm{u}+\big(\rho(\phi)\bm{u}-\rho_dm_\epsilon(\phi)\nabla\mu\big)\cdot\nabla\bm{u} & =\nabla\cdot\big(2\eta(\phi)\D(\bm{u})-p\I\big)+\mu\nabla\phi+\bm{J}\times\bm{B}+\rho\bm{g}, \label{phase_model_3} \\
        \nabla\cdot\bm{u} & =0, \label{phase_model_4} \\
        \bm{J} & =\sigma(\phi)\big(\bm{E}+\bm{u}\times\bm{B}\big), \label{phase_model_5} \\
        \bm{E} & = -\nabla{V}, \label{tmp_E}\\
        \nabla\cdot\bm{J} & =0. \label{phase_model_6}
    \end{align}
\end{subequations}
In the above system, $\phi$ is the order parameter that labels the two fluids such that
\begin{equation*}
    \phi(\bm{x},t)=
    \begin{cases}
        1,  & \bm{x}\in\Omega_+(t),\\
        -1, & \bm{x}\in\Omega_-(t).
    \end{cases}
\end{equation*}
Additionally, $\mu$ is the chemical potential serving as an auxiliary variable, $m_\epsilon(\phi)>0$ represents the mobility which plays the role of diffusivity for the diffuse interface, and $\widehat{\lambda}=\frac{3\lambda}{2\sqrt{2}}$ stands for the scaled surface tension. Moreover, the density $\rho$, the dynamic viscosity $\eta$, and the electrical conductivity $\sigma$ of the mixture are all linearly dependent on $\phi$ such that (exemplarily for $\rho$)
\begin{equation*}
    \rho=\frac{\rho_+-\rho_-}{2}\phi+\frac{\rho_++\rho_-}{2}.
\end{equation*}
We define $\rho_d:=\partial\rho/\partial\phi=(\rho_+-\rho_-)/2$ as the density difference between the two fluids. Equations \eqref{phase_model_1}-\eqref{phase_model_2} represent the widely recognized Cahn-Hilliard equation which describes the process of phase separation. It is necessary to point out that equation \eqref{tmp_E} represents merely a constitutive relation, and $\bm{E}$ will not be used in both the asymptotic analysis and the numerical computations.

To finalize the system, we impose some appropriate initial and boundary conditions. The initial data is set as 
\begin{equation*}
    \bm{u}(0)=\bm{u}^0	
    \qquad	\mbox{and} \qquad
     \phi(0)=\phi^0.
\end{equation*}
We only consider the enclosed flows in this context, and therefore, no-slip walls are imposed for $\bm{u}$, yielding
\begin{equation*}
    \bm{u}=\bm{0} \qquad \text{on}\;\;\partial\Omega,
\end{equation*}
and insulating walls are applied for $\bm{J}$ and $V$, leading to
\begin{equation*}
    \bm{J}\cdot\bm{n}=0 \qquad \text{on}\;\;\partial\Omega,
\end{equation*}
where $\bm{n}$ is the unit outward normal to $\Omega$. Mass conservation gives no-flux walls for $\phi$ and $\mu$
\begin{equation*}
    \begin{aligned}
        \nabla\phi\cdot\bm{n} & =0 \qquad \text{on}\;\;\partial\Omega, \\
        \nabla\mu\cdot\bm{n} & =0 \qquad \text{on}\;\;\partial\Omega.
    \end{aligned}
\end{equation*}

\begin{remark}
\label{energy_relation}
    It can be readily confirmed that the sharp interface model \eqref{sharp_interface} satisfies the individual mass conservation and the energy dissipation law
    \begin{equation*}
        \begin{aligned}
            \frac{\mathrm{d}}{\mathrm{d}t}\int_{\Omega_\pm(t)}\rho_\pm\mathrm{d}\bm{x} & =0, \\
            \frac{\mathrm{d}}{\mathrm{d}t}\left(\int_\Omega\frac{1}{2}\rho|\bm{u}|^2\mathrm{d}\bm{x}+\int_{\Gamma(t)}\lambda\mathrm{d}\bm{s}+\int_\Omega\rho{g}z\mathrm{d}\bm{x}\right) & =-\int_\Omega{2}\eta|\D(\bm{u})|^2\mathrm{d}\bm{x}-\int_\Omega\sigma^{-1}|\bm{J}|^2\mathrm{d}\bm{x}\le{0},
        \end{aligned}
    \end{equation*}
    with $\bm{x}=(x, y, z)^{\top}$ and $\rho=\rho_+\mathscr{X}_{\Omega_+}+\rho_-\mathscr{X}_{\Omega_-}$. Nevertheless, for the proposed diffuse interface model \eqref{phase_field}, only the global mass conservation and the subsequent energy identity are valid,
    \begin{equation*}
        \begin{aligned}
            \frac{\mathrm{d}}{\mathrm{d}t}\int_\Omega\phi\mathrm{d}\bm{x} & =0, \\
            \frac{\mathrm{d}}{\mathrm{d}t}\left(\int_\Omega\frac{1}{2}\rho|\bm{u}|^2\mathrm{d}\bm{x}+\int_\Omega\frac{\widehat{\lambda}\epsilon}{2}|\nabla\phi|^2\mathrm{d}\bm{x}+\int_\Omega\frac{\widehat{\lambda}\epsilon^{-1}}{4}\left(\phi^2-1\right)^2\mathrm{d}\bm{x}\right) & =-\int_\Omega m_\epsilon|\nabla\mu|^2\mathrm{d}\bm{x}-\int_\Omega{2}\eta|\D(\bm{u})|^2\mathrm{d}\bm{x} \\
            & \;\;\;\;-\int_\Omega\sigma^{-1}|\bm{J}|^2\mathrm{d}\bm{x}+\int_\Omega\rho\bm{g}\cdot\bm{u}\mathrm{d}\bm{x}.
        \end{aligned}
    \end{equation*}
    The practical applications usually involve the gravity, and the energy of the diffuse interface model \eqref{phase_field} is no longer dissipative as that of the sharp interface counterpart. This unexpected energy law fundamentally attributes to the global incompressibility constraint \eqref{phase_model_4}. Simultaneously, in order to guarantee the thermodynamic consistency at a theoretical level, an extra flux $\rho_dm_\epsilon\nabla\mu$ must be incorporated into the momentum equation \eqref{phase_model_3}. In turn, this embarrasses the convergence analysis and the design of efficient structure-preserving schemes (a numerical demonstration on the influence of this flux can be found in \cite{2014_Grun}). However, the proposed model is capable of simulating MHD effects in general engineering processes, and remains a solenoidal velocity as well. In a contrast, interested readers are referred to \cite{1998_Lowengrub} for another modeling approach leading to a much more complex quasi-incompressible diffuse interface model.
\end{remark}

It is highly desirable to ensure the correct physical behavior for the simulation of diffuse interface models, compared to that of the sharp interface system \eqref{sharp_interface}. In other words, whether the standard sharp interface formulation \eqref{sharp_interface} can be recovered from the diffuse interface model \eqref{phase_field} as the interface thickness $\epsilon$ approaches zero becomes a critical issue. In addition, it is known that the scaling between $\epsilon$ and the mobility $m_\epsilon$ is also critical during the simulations. In this study, we take a constant mobility with a scaling law $m_\epsilon=\mathcal{O}(\epsilon^2)$ employed in \cite{2013_Magaletti, 2023_Khanwale, 2024_Wang}. This scaling law facilitates the asymptotic analysis by making the extra flux term $m_\epsilon\nabla\mu$ in the momentum equation a higher-order contribution with respect to $\epsilon$. Consequently, this term does not affect the leading-order dynamics, ensuring that the correct physical behavior is preserved in the sharp interface asymptotics.

\begin{theorem}
    \label{sharp_asymptotics}
    Assuming the scaling law $m_\epsilon=\mathcal{O}(\epsilon^2)$, the sharp interface limit of model \eqref{phase_field} is system \eqref{sharp_interface}.
\end{theorem}

The derivation of the diffuse interface model \eqref{phase_field}, and the corresponding asymptotic analysis on the relationship between the system \eqref{phase_field}  and \eqref{sharp_interface} will be given in the following two subsections, respectively.

\subsection{Model derivation}
\label{model_derivation}

In this subsection, we present the derivation of the system \eqref{phase_field} by means of the Onsager's variational principle and the laws of thermodynamics. First, we assume a partial mixing of the two fluids in a narrow interfacial zone of thickness $\epsilon$ with the incompressibility on the volume. Let $V(t)\subset\Omega$ be an arbitrary control volume with its mass and density being given by $M$ and $\rho$. Then, denoted by $M_\pm$ and $\varrho_\pm$ the masses and densities of the two fluids in $V$, it holds that
 \begin{equation*}
    \varrho_\pm=\frac{M_\pm}{V}
    \qquad \mbox{and} \qquad
    \rho=\frac{M}{V}=\frac{M_++M_-}{V}=\varrho_++\varrho_-.
\end{equation*}
In turn, the volume incompressibility assumption leads to
\begin{equation*}
    1=\frac{\frac{M_+}{\rho_+}+\frac{M_-}{\rho_-}}{V}=\frac{\varrho_+}{\rho_+}+\frac{\varrho_-}{\rho_-}.
\end{equation*}
We define the order parameter for the mixture as
\begin{equation*}
    \phi=\frac{\varrho_+}{\rho_+}-\frac{\varrho_-}{\rho_-},
\end{equation*}
which indicates that
\begin{equation*}
    \rho=\frac{\rho_+-\rho_-}{2}\phi+\frac{\rho_++\rho_-}{2}.
\end{equation*}
Supposing that the fluids move with different velocities $\bm{u}_\pm$, the continuity equations in the bulk read
\begin{equation*}
    \partial_t\varrho_\pm+\nabla\cdot\left(\varrho_\pm\bm{u}_\pm\right)=0.
\end{equation*}
Then we introduce the volume-averaged velocity and relative mass flux as
\begin{equation*}
    \bm{u}=\frac{\varrho_+}{\rho_+}\bm{u}_++\frac{\varrho_-}{\rho_-}\bm{u}_-
    \qquad \mbox{and} \qquad
    \bm{j}_\pm=\varrho_\pm\left(\bm{u}_\pm-\bm{u}\right).
\end{equation*}
As a consequence, we obtain
\begin{equation*}
    \nabla\cdot\bm{u}=\nabla\cdot\left(\frac{\varrho_+}{\rho_+}\bm{u}_++\frac{\varrho_-}{\rho_-}\bm{u}_-\right)=-\partial_t\left(\frac{\varrho_+}{\rho_+}+\frac{\varrho_-}{\rho_-}\right)=0,
\end{equation*}
and
\begin{equation*}
    \partial_t\phi=-\nabla\cdot\left(\frac{\varrho_+}{\rho_+}\bm{u}+\frac{\bm{j}_+}{\rho_+}-\frac{\varrho_-}{\rho_-}\bm{u}-\frac{\bm{j}_-}{\rho_-}\right)=-\nabla\cdot\left(\phi\bm{u}+\bm{j}_\phi\right),
\end{equation*}
where $\bm{j}_\phi=\frac{\bm{j}_+}{\rho_+}-\frac{\bm{j}_-}{\rho_-}$. Correspondingly, this implies that the mass diffusion is permitted in our method. Hence, the momentum $\rho\bm{u}$ should be transported by $\rho\bm{u}+\rho_d\bm{j}_\phi$, giving
\begin{equation*}
    \rho\partial_t\bm{u}+\left(\rho\bm{u}+\rho_d\bm{j}_\phi\right)\cdot\nabla\bm{u}=\nabla\cdot\S-\nabla{p}+\bm{F},
\end{equation*}
where $\S$ is a symmetric second-order tensor, $p$ is the pressure serving as a Lagrange multiplier for the incompressibility constraint, and $\bm{F}$ is a force density. Invoking Ohm's law with the charge conservation, we conclude
\begin{subequations}
    \begin{align*}
        \partial_t\phi+\nabla\cdot\left(\phi\bm{u}+\bm{j}_\phi\right) & =0, \\
        \rho\partial_t\bm{u}+\left(\rho\bm{u}+\rho_d\bm{j}_\phi\right)\cdot\nabla\bm{u} & =\nabla\cdot\S-\nabla{p}+\bm{F}, \\
        \nabla\cdot\bm{u} & =0, \\
        \bm{J} & =\sigma(\phi)\left(\bm{E}+\bm{u}\times\bm{B}\right), \\
        \bm{E} & =-\nabla{V}, \\
        \nabla\cdot\bm{J} & =0,
    \end{align*}
\end{subequations}
where $\bm{j}_\phi$, $\S$, and $\bm{F}$ are to be determined, and the definition of $\sigma(\phi)$ will be given later.

To derive a thermodynamically consistent model, we employ the Onsager's variational principle as in \cite{2012_Abels}. By introducing the Rayleighian functional
\begin{equation*}
    \mathcal{R}=\frac{\mathrm{d}\mathcal{E}}{\mathrm{d}t}+\mathcal{P},
\end{equation*}
where $\mathcal{E}$ and $\mathcal{P}$ are the free energy and dissipation functional, the Onsager's variational principle states that the kinetic equation is equivalent to the minimum of $\mathcal{R}$ \cite{2006_Qian}. Then, we define
\begin{equation*}
    \mathcal{E}=\int_\Omega\left(\frac{1}{2}\rho|\bm{u}|^2+f(\phi,\nabla\phi)\right)\mathrm{d}\bm{x},
\end{equation*}
where $f(\phi,\nabla\phi)$ is a Helmholtz free energy density. Thus, it holds that
\begin{equation*}
    \frac{\mathrm{d}\mathcal{E}}{\mathrm{d}t}=\int_\Omega\left(-\S:\D(\bm{u})+\bm{u}\cdot\bm{F}-\mu\bm{u}\cdot\nabla\phi+\bm{j}_\phi\cdot\nabla\mu\right)\mathrm{d}\bm{x},
\end{equation*}
where $\mu=\frac{\partial{f}}{\partial\phi}-\nabla\cdot\frac{\partial{f}}{\partial\nabla\phi}$ is abbreviated to denote the chemical potential, and here we assume a vanishing normal component of $\frac{\partial{f}}{\partial\nabla\phi}$ on $\partial\Omega$. Since the first law of thermodynamics reads
\begin{equation*}
    \frac{\mathrm{d}\mathcal{E}}{\mathrm{d}t}=-\mathcal{T}\frac{\mathrm{d}\mathcal{S}}{\mathrm{d}t}+\frac{\mathrm{d}\mathcal{W}}{\mathrm{d}t},
\end{equation*}
where $\mathcal{T}$, $\mathcal{S}$, and $\mathcal{W}$ represent the temperature, entropy, and mechanical work, we are able to identify
\begin{equation*}
    \frac{\mathrm{d}\mathcal{W}}{\mathrm{d}t}=\int_\Omega\left(\bm{u}\cdot\bm{F}-\mu\bm{u}\cdot\nabla\phi\right)\mathrm{d}\bm{x}.
\end{equation*}
Incorporating external gravitational and Lorentz forces, we can write $\bm{F}$ as
\begin{equation*}
    \bm{F}=\mu\nabla\phi+\bm{J}\times\bm{B}+\rho\bm{g}.
\end{equation*}
To further specify $\S$ and $\bm{j}_\phi$, we introduce the dissipation functional
\begin{equation*}
    \mathcal{P}=\int_\Omega\left(\frac{|\S|^2}{4\eta(\phi)}+\frac{|\bm{j}_\phi|^2}{2m_\epsilon(\phi)}+\frac{|\bm{J}|^2}{2\sigma(\phi)}\right)\mathrm{d}\bm{x},
\end{equation*}
which accounts for the viscous dissipation, the diffusive dissipation, and the Ohmic dissipation, respectively. In accordance to the Onsager's variational principle, we get
\begin{equation*}
    \S=2\eta\D(\bm{u})
    \qquad \mbox{and} \qquad
    \bm{j}_\phi=-m_\epsilon\nabla{\mu}.
\end{equation*}
In this work, we consider the Helmholtz free energy density in the Ginzburg-Landau form:
\begin{equation*}
    f(\phi,\nabla\phi)=\frac{\widehat{\lambda}\epsilon}{2}|\nabla\phi|^2+\frac{\widehat{\lambda}\epsilon^{-1}}{4}\left(\phi^2-1\right)^2.
\end{equation*}
Moreover, we assume that the material properties vary linearly in the mixture, i.e.,
\begin{equation*}
    \eta=\frac{\eta_+-\eta_-}{2}\phi+\frac{\eta_++\eta_-}{2}
    \qquad \mbox{and} \qquad
    \sigma=\frac{\sigma_+-\sigma_-}{2}\phi+\frac{\sigma_++\sigma_-}{2},
\end{equation*}
so that the model is compatible with the inductionless MHD equations in the single fluid regime. This finally leads to the system \eqref{phase_field}.

\subsection{The proof of Theorem \ref{sharp_asymptotics}}

We exploit the method of formally matched asymptotic expansions to identify the sharp interface limit of system \eqref{phase_field}. This argument depends on two critical assumptions: (i) with vanishing interface thickness the domain can be decomposed into two distinct subdomains separated by an interface, each containing only one fluid; and (ii) the solutions admit separate asymptotic expansions in terms of the interface thickness in the bulk and near the interface, which must be matched in the overlapping region.

\subsubsection{Preliminaries}

For a solution $\left(\phi_\epsilon,\mu_\epsilon,\bm{u}_\epsilon,p_\epsilon,\bm{J}_\epsilon,V_\epsilon\right)$ of system \eqref{phase_field}, we assume that it converges formally to a limit $\left(\phi,\mu,\bm{u},p,\bm{J},V\right)$ as the interface thickness $\epsilon\to{0}$, and further,
\begin{equation*}
    \Omega_+(t)=\left\{\bm{x}\in\Omega:\phi_\epsilon(\bm{x})>0\right\}
    \qquad \mbox{and} \qquad
    \Omega_-(t)=\left\{\bm{x}\in\Omega:\phi_\epsilon(\bm{x})<0\right\}.
\end{equation*}
Then, the corresponding interface separating $\Omega_\pm(t)$ is denoted by
\begin{equation*}
    \Gamma(t)=\left\{\bm{x}\in\Omega:\phi_\epsilon(\bm{x})=0\right\}.
\end{equation*}

In outer regions far from $\Gamma(t)$, there holds the outer expansions with respect to the interface thickness $\epsilon$ of $\phi_\epsilon$, $\bm{u}_\epsilon$, $p_\epsilon$, $\bm{J}_\epsilon$, $V_\epsilon$ in the form (exemplarily for $\phi_\epsilon$)
\begin{equation*}
    \phi_\epsilon=\phi_0^\pm+\epsilon\phi_1^\pm+\epsilon^2\phi_2^\pm+\cdots.
\end{equation*}
As an auxiliary variable, we have evidently for $\mu$:
\begin{equation*}
    \mu_\epsilon=\epsilon^{-1}\mu_0^\pm+\mu_1^\pm+\epsilon\mu_2^\pm+\cdots
\end{equation*}
with $\mu_0^\pm=\widehat{\lambda}\left((\phi_0^\pm)^3-\phi_0^\pm\right)$. In the transition layer near $\Gamma(t)$, we introduce a new stretched coordinates to scale the variables. Let $d(\bm{x})$ be the signed distance function to $\Gamma(t)$ such that $d(\bm{x})>0$ if $\bm{x}\in\Omega_+$, then we have $\bm{n}_\Gamma=\nabla{d}$ and $V_\Gamma=-\partial_td$. Defining the scaled distance $z=d/\epsilon$, the following change of variables holds \cite{2012_Abels}:
\begin{subequations}
\label{change_variable}
    \begin{align}
        \partial_t\chi & =-\epsilon^{-1}V_\Gamma\partial_z\chi+\mathcal{O}(1), \\
        \nabla\chi & =\epsilon^{-1}\partial_z\chi\bm{n}_\Gamma+\nabla_\Gamma\chi+\mathcal{O}(\epsilon), \\
        \nabla\cdot\bm{\chi} & =\epsilon^{-1}\partial_z\bm{\chi}\cdot\bm{n}_\Gamma+\nabla_\Gamma\cdot\bm{\chi}+\mathcal{O}(\epsilon), \\
        \Delta\chi & =\epsilon^{-2}\partial_{zz}\chi-\epsilon^{-1}\kappa\partial_z\chi+\mathcal{O}(1),
    \end{align}
\end{subequations}
where $\chi$ and $\bm{\chi}$ represent the generic scalar-valued and vector-valued variables, respectively. Then in this region we have the inner expansions in terms of the interface thickness $\epsilon$ for $\phi_\epsilon$, $\bm{u}_\epsilon$, $p_\epsilon$, $\bm{J}_\epsilon$, $V_\epsilon$ in the form (also exemplarily for $\phi_\epsilon$)
\begin{equation*}
    \phi_\epsilon=\phi_0^i+\epsilon\phi_1^i+\epsilon^2\phi_2^i+\cdots,
\end{equation*}
with $\mu_\epsilon=\epsilon^{-1}\mu_0^i+\mu_1^i+\epsilon\mu_2^i+\cdots$, where
\begin{equation*}
    \mu_0^i=\widehat{\lambda}\left(-\partial_{zz}\phi_0^i+(\phi_0^i)^3-\phi_0^i\right)
    \qquad \mbox{and} \qquad
    \mu_1^i=\widehat{\lambda}\left(-\partial_{zz}\phi_1^i+\kappa\partial_z\phi_0^i+(3(\phi_0^i)^2-1)\phi_1^i\right).
\end{equation*}
In the overlapping region, the outer and inner expansions must satisfy the matching conditions \cite{2012_Abels}
\begin{subequations}
    \begin{align}
        \lim_{z\to\pm\infty}\chi_0^i & =\chi_0^\pm,\label{match_1} \\
        \lim_{z\to\pm\infty}\partial_z\chi_0^i & =0,\label{match_2} \\
        \lim_{z\to\pm\infty}\partial_z\chi_1^i & =\nabla\chi_0^\pm\cdot\bm{n}_\Gamma.\label{match_3}
    \end{align}
\end{subequations}

Next, we will establish that the limit solution $\left(\phi,\mu,\bm{u},p,\bm{J},V\right)$ exactly satisfies the sharp interface model \eqref{sharp_interface} by following the procedure in \cite{2012_Abels, 2018_Xu}.

\subsubsection{Bulk equations}

\textbf{Expansion of \eqref{phase_model_3} at order $\mathcal{O}(\epsilon^{-1})$:} Recalling the assumption $m_\epsilon:=m\epsilon^2$, we naturally obtain that the leading order of the flux $m_\epsilon\nabla\mu$ is $\mathcal{O}(\epsilon)$. Hence, we arrive at
\begin{equation*}
    \mu_0^\pm\nabla\phi_0^\pm=0,
\end{equation*}
which is equivalent to
\begin{equation*}
    \nabla\frac{\left((\phi_0^\pm)^2-1\right)^2}{4}=0.
\end{equation*}
This implies that $\phi_0^\pm=c_\pm$ with $c_\pm$ two constants such that $c_+>0$ and $c_-<0$.

\textbf{Expansion of \eqref{phase_model_3}-\eqref{phase_model_6} at order $\mathcal{O}(1)$:} Using the fact that $m_\epsilon\nabla\mu$ is of order $\mathcal{O}(\epsilon)$, we conclude
\begin{equation}
\label{bulk_equation}
    \begin{aligned}
        \rho_0^\pm\left(\partial_t\bm{u}_0^\pm+\bm{u}_0^\pm\cdot\nabla\bm{u}_0^\pm\right) & =\nabla\cdot\big(2\eta_0^\pm\D(\bm{u}_0^\pm)-\nabla{p}_0^\pm\I\big)+\mu_0^\pm\nabla\phi_1^\pm+\bm{J}_0^\pm\times\bm{B}+\rho_0^\pm\bm{g}, \\
        \nabla\cdot\bm{u}_0^\pm & =0, \\
        \bm{J}_0^\pm & =\sigma_0^\pm\left(\bm{E}^\pm+\bm{u}_0^\pm\times\bm{B}\right), \\
        \bm{E}^\pm &= -\nabla{V}_0^\pm, \\
        \nabla\cdot\bm{J}_0^\pm & = 0,
    \end{aligned}
\end{equation}
where we abbreviate $\chi_0^\pm=\frac{\chi_+-\chi_-}{2}\phi_0^\pm+\frac{\chi_++\chi_-}{2}$ for $\chi\in\{\rho, \eta, \sigma\}$.

\subsubsection{Jump conditions}

\textbf{Expansion of \eqref{phase_model_4} at order $\mathcal{O}(\epsilon^{-1})$:} Invoking change of variables \eqref{change_variable}, we get
\begin{equation}
\label{inner_1}
    \partial_z\bm{u}_0^i\cdot\bm{n}_\Gamma=0,
\end{equation}
which implies that $\bm{u}_0^i\cdot\bm{n}_\Gamma$ is independent of $z$.

\textbf{Expansion of \eqref{phase_model_1} at order $\mathcal{O}(\epsilon^{-1})$:} With a similar procedure, we obtain
\begin{equation*}
    -V_\Gamma\partial_z\phi_0^i+\partial_z\left(\phi_0^i\bm{u}_0^i\right)\cdot\bm{n}_\Gamma=m\partial_{zz}\mu_0^i.
\end{equation*}
Integrating this equation with respect to $z$ from $-\infty$ to $\infty$ yields
\begin{equation*}
    (-V_\Gamma+\bm{u}_0^i\cdot\bm{n}_\Gamma)\int_{-\infty}^\infty\partial_z\phi_0^i\mathrm{d}z=m\int_{-\infty}^\infty\partial_{zz}\mu_0^i\mathrm{d}z.
\end{equation*}
Then matching conditions \eqref{match_1}-\eqref{match_2} requires
\begin{equation*}
    \left(-V_\Gamma+\bm{u}_0^i\cdot\bm{n}_\Gamma\right)\left(\phi_0^+-\phi_0^-\right)=0,
\end{equation*}
in turn giving
\begin{equation}
    V_\Gamma=\bm{u}_0^i\cdot\bm{n}_\Gamma.
\end{equation}
Additionally, we have
\begin{equation*}
    \partial_{zz}\mu_0^i=0,
\end{equation*}
or equivalently, $\mu_0^i$ linearly depends on $z$. Required by matching condition \eqref{match_1}, we conclude that $\mu_0^i$ has to be a constant.

\textbf{Expansion of \eqref{phase_model_3} at order $\mathcal{O}(\epsilon^{-2})$:} Denoting by $\mathcal{D}(\bm{A})=\frac{1}{2}\left(\bm{A}+\bm{A}^\top\right)$, we obtain after the change of variables \eqref{change_variable},
\begin{equation}
    \nabla\cdot\left(2\eta\D(\bm{u})\right)=\epsilon^{-2}\partial_z\left(2\eta\mathcal{D}(\partial_z\bm{u}\otimes\bm{n}_\Gamma)\bm{n}_\Gamma\right)+\epsilon^{-1}\partial_z\left(2\eta\mathcal{D}(\nabla_\Gamma\bm{u})\bm{n}_\Gamma\right)+\epsilon^{-1}\nabla_\Gamma\cdot\left(2\eta\mathcal{D}(\partial_z\bm{u}\otimes\bm{n}_\Gamma\right)+\mathcal{O}(1).
\end{equation}
Notice that $m_\epsilon\nabla\mu$ is of order $\mathcal{O}(\epsilon)$ in the inner region since $\mu_0^i$ is independent of $z$, we conclude
\begin{equation*}
    \partial_z\left(\eta_0^i\partial_z\bm{u}_0^i\right)+\mu_0^i\partial_z\phi_0^i\bm{n}_\Gamma=0,
\end{equation*}
where we have used
\begin{equation*}
    \left(\bm{n}_\Gamma\otimes\partial_z\bm{u}_0^i\right)\bm{n}_\Gamma=\left(\partial_z\bm{u}_0^i\cdot\bm{n}_\Gamma\right)\bm{n}_\Gamma=0.
\end{equation*}
By integrating this equation with respect to $z$, we can obtain
\begin{equation*}
    \left.\eta_0^i\partial_z\bm{u}_0^i\right|_{-\infty}^\infty+\mu_0^i\left.\phi_0^i\right|_{-\infty}^\infty\bm{n}_\Gamma=0.
\end{equation*}
Through equation \eqref{inner_1}, taking an inner product with $\bm{n}_\Gamma$ and matching condition \eqref{match_1} yield
\begin{equation*}
    \mu_0^i\left(\phi_0^+-\phi_0^-\right)=0,
\end{equation*}
in turn indicating $\mu_0^i=0$, or equivalently
\begin{equation*}
    -\partial_{zz}\phi_0^i+\left(\phi_0^i\right)^3-\phi_0^i=0.
\end{equation*}
The solvability of this equation \cite{2018_Xu} leads to
\begin{equation*}
    \phi_0^\pm=\lim_{z\to\pm\infty}\phi_0^i=\pm{1}\;\;\text{and }\;\;\phi_0^i=\tanh\frac{z}{\sqrt{2}}.
\end{equation*}
Then we have a second-order ordinary differential equation
\begin{equation*}
    \partial_z\left(\eta_0^i\partial_z\bm{u}_0^i\right)=0,
\end{equation*}
which, by matching conditions \eqref{match_1}-\eqref{match_2} and positivity of $\eta$, implies that $\bm{u}_0^i$ is independent of $z$. Hence, we obtain
\begin{equation*}
    [\![\bm{u}_0^\pm]\!]=0.
\end{equation*}

\textbf{Expansion of \eqref{phase_model_3} at order $\mathcal{O}(\epsilon^{-1})$:} The combination of matching conditions \eqref{match_1} and \eqref{match_3} yields 
\begin{equation}
\label{match_4}
    \lim_{z\to\pm\infty}\left(\partial_z\bm{u}_1^i\otimes\bm{n}_\Gamma+\nabla_\Gamma\bm{u}_0^i\right)=\nabla\bm{u}_0^\pm.
\end{equation}
Notice that $\bm{u}_0^i$ and $\mu_0^i$ are indeed constants, we can derive
\begin{equation*}
    \partial_z\left(2\eta_0^i\mathcal{D}(\partial_z\bm{u}_1^i\otimes\bm{n}_\Gamma)\bm{n}_\Gamma\right)+\partial_z\left(2\eta_0^i\mathcal{D}(\nabla_\Gamma\bm{u}_0^i)\bm{n}_\Gamma\right)-\partial_zp_0^i\bm{n}_\Gamma+\mu_1^i\partial_z\phi_0^i\bm{n}_\Gamma=0.
\end{equation*}
Then integrating and using matching conditions \eqref{match_1} and \eqref{match_4} yield
\begin{equation*}
    [\![2\eta_0^\pm\D(\bm{u}_0^\pm)-p_0^\pm\I]\!]\bm{n}_\Gamma=-\int_{-\infty}^\infty\mu_1^i\partial_z\phi_0^i\bm{n}_\Gamma\mathrm{d}z.
\end{equation*}
For the right-hand side of this equation, integration by parts with matching conditions \eqref{match_1}-\eqref{match_2} leads to
\begin{equation*}
    \begin{aligned}
        \int_{-\infty}^\infty\mu_1^i\partial_z\phi_0^i\mathrm{d}z & =\widehat{\lambda}\int_{-\infty}^\infty\left(-\partial_{zz}\phi_1^i+\kappa\partial_z\phi_0^i+(3(\phi_0^i)^2-1)\phi_1^i\right)\partial_z\phi_0^i\mathrm{d}z, \\
        & =\widehat{\lambda}\kappa\int_{-\infty}^\infty(\partial_z\phi_0^i)^2\mathrm{d}z+\widehat{\lambda}\int_{-\infty}^\infty\left(-\partial_{zz}\phi_1^i+(3(\phi_0^i)^2-1)\phi_1^i\right)\partial_z\phi_0^i\mathrm{d}z \\
        & =\frac{2\sqrt{2}\widehat{\lambda}\kappa}{3}+\widehat{\lambda}\int_{-\infty}^\infty\left(\partial_{zz}\phi_0^i-(\phi_0^i)^3+\phi_0^i\right)\partial_{z}\phi_1^i\mathrm{d}z\\
        & =\lambda\kappa,
    \end{aligned}
\end{equation*}
where
\begin{equation*}
    \int_{-\infty}^\infty(\partial_z\phi_0^i)^2\mathrm{d}z=\frac{1}{2}\int_{-\infty}^\infty\mathrm{sech}^4\frac{z}{\sqrt{2}}\mathrm{d}z=\frac{2\sqrt{2}}{3}.
\end{equation*}
Then we conclude
\begin{equation*}
    [\![2\eta_0^\pm\D(\bm{u}_0^\pm)-p_0^\pm\I]\!]\bm{n}_\Gamma=-\lambda\kappa.
\end{equation*}

\textbf{Expansion of \eqref{phase_model_5}-\eqref{phase_model_6} at order $\mathcal{O}(\epsilon^{-1})$:} A direct change of variables \eqref{change_variable} gives
\begin{equation*}
    \begin{aligned}
        \partial_z\varphi_0^i\bm{n}_\Gamma & =\bm{0}, \\
        \partial_z\bm{J}_0^i\cdot\bm{n}_\Gamma & =0.
    \end{aligned}
\end{equation*}
With matching condition \eqref{match_1}, we obtain
\begin{subequations}
    \begin{align*}
        [\![V_0^\pm]\!] & =0, \\
        [\![\bm{J}_0^\pm\cdot\bm{n}_\Gamma]\!] & =0.
    \end{align*}
\end{subequations}

\subsubsection{Summary of the sharp interface limit}

We summarize the main results as follows. For the region far from $\Gamma(t)$, we have $\phi_0^\pm=\pm{1}$ and $\mu_0^\pm=0$. Then, the bulk equations \eqref{bulk_equation} can be simplified as
\begin{subequations}
    \begin{align*}
        \rho_\pm\left(\partial_t\bm{u}_0^\pm+\bm{u}_0^\pm\cdot\nabla\bm{u}_0^\pm\right) & =\nabla\cdot\left(2\eta_0^\pm\D(\bm{u}_0^\pm)-p_0^\pm\I\right)+\bm{J}_0^\pm\times\bm{B}+\rho_0^\pm\bm{g}, \\
        \nabla\cdot\bm{u}_0^\pm & =0, \\
        \bm{J}_0^\pm & =\sigma_\pm\left(\bm{E}^\pm+\bm{u}_0^\pm\times\bm{B}\right), \\
        \bm{E}^\pm & =-\nabla{V}_0^\pm, \\
        \nabla\cdot\bm{J}_0^\pm & = 0,
    \end{align*}
\end{subequations}
and the jump conditions are
\begin{subequations}
    \begin{align*}
        [\![2\eta_0^\pm\D(\bm{u}_0^\pm)-p_0^\pm\I]\!]\bm{n}_\Gamma & =-\lambda\kappa, \\
        [\![\bm{u}_0^\pm]\!] & =\bm{0}, \\
        [\![\bm{J}_0^\pm\cdot\bm{n}_\Gamma]\!] & =0, \\
        [\![V_0^\pm]\!] & =0, \\
        V_\Gamma & =\bm{u}_0^i\cdot\bm{n}_\Gamma.
    \end{align*}
\end{subequations}
We can deduce that the limit system is the sharp interface model \eqref{sharp_interface} as the interface thickness $\epsilon$ approaches zero.

\section{Numerical results}
\label{section_3}

\subsection{Non-dimensionalization and discrete scheme}

Let $L_r$, $u_r$, $\rho_r$, $\eta_r$, $\sigma_r$ denote the characteristic quantities of length, velocity, density, dynamic viscosity, and electrical conductivity, respectively. We introduce the following scalings:
\begin{equation*}
    \begin{aligned}
         & \bm{x}\leftarrow\frac{\bm{x}}{L_r},\;\;{t}\leftarrow\frac{tu_r}{L_r},\;\;\bm{u}\leftarrow\frac{\bm{u}}{u_r},\;\;\rho\leftarrow\frac{\rho}{\rho_r},\;\;\eta\leftarrow\frac{\eta}{\eta_r},\;\;{p}\leftarrow\frac{p}{\rho_ru_r^2}, \\
         & \mu\leftarrow\frac{\mu\epsilon}{\widehat{\lambda}},\;\;{V}\leftarrow\frac{{V}}{L_ru_rB_r},\;\;\bm{J}\leftarrow\frac{\bm{J}}{\sigma_ru_rB_r},\;\;\sigma\leftarrow\frac{\sigma}{\sigma_r},\;\;\bm{E}\leftarrow\frac{\bm{E}}{u_rB_r}.
    \end{aligned}
\end{equation*}
Certainly, the dimensionless system is given by
\begin{subequations}
    \begin{align*}
        \partial_t\phi+\nabla\cdot\left(\phi\bm{u}\right) & =\frac{1}{\mathrm{Pe}}\Delta\mu, \\
        \mu & =-\mathrm{Cn}^2\Delta\phi+\phi^3-\phi, \\
        \rho\partial_t\bm{u}+\left(\rho\bm{u}-\rho_d\frac{1}{\mathrm{Pe}}\nabla\mu\right)\cdot\nabla\bm{u} & =\frac{1}{\mathrm{Re}}\nabla\cdot\left(2\eta\D(\bm{u})\right)-\nabla{p}+\frac{1}{\mathrm{We}\mathrm{Cn}}\mu\nabla\phi+\mathrm{N}\bm{J}\times\bm{B}-\frac{1}{\mathrm{Fr}}\rho\bm{e}_3, \\
        \nabla\cdot\bm{u} & =0, \\
        \bm{J} & =\sigma\left(\bm{E}+\bm{u}\times\bm{B}\right), \\
        \bm{E} & =-\nabla{V}, \\
        \nabla\cdot\bm{J} & =0,
    \end{align*}
\end{subequations}
where $\bm{e}_3=(0,0,1)^{\top}$. Here, the dimensionless numbers are
\begin{itemize}
    \item Cahn number $\mathrm{Cn}=\displaystyle\frac{\epsilon}{L_r}$ (interface thickness);
    \item P\'{e}clet number $\mathrm{Pe}=\displaystyle\frac{\epsilon{L}_ru_r}{\widehat{\lambda}m_\epsilon}$ (ratio of interface thickness to mobility);
    \item Reynolds number $\mathrm{Re}=\displaystyle\frac{L_r\rho_ru_r}{\eta_r}$ (ratio of the inertial to viscous forces);
    \item Weber number $\mathrm{We}=\displaystyle\frac{L_r\rho_ru_r^2}{\widehat{\lambda}}$ (ratio of inertial to interfacial forces);
    \item Stuart number $\mathrm{N}=\displaystyle\frac{L_r\sigma_rB_r^2}{\rho_ru_r}$ (ratio of electromagnetic to inertial forces);
    \item Froude number $\mathrm{Fr}=\displaystyle\frac{u_r^2}{gL_r}$ (ratio of inertial to gravitational acceleration forces).
\end{itemize}
In this scaling, we apply the scaling law $\frac{1}{\mathrm{Pe}}=3\mathrm{Cn}$ given by Magaletti et al. \cite{2013_Magaletti}.

Let $L^2(\Omega)$ be the Lebesgue space of square integrable functions equipped with inner product $(\cdot,\cdot)$, and let $H^1(\Omega)$ and $\bm{H}(\mathrm{div};\Omega)$ be its subspaces with square integrable gradients and square integrable divergences, respectively. We define their subspaces $L_0^2(\Omega)$, $H_0^1(\Omega)$, and $\bm{H}_0(\mathrm{div};\Omega)$ with vanishing mean values on $\Omega$, vanishing traces, and vanishing normal traces on $\partial\Omega$. Then let $Q_h\subset{H}^1(\Omega)$, $\bm{X}_h\subset{H}_0^1(\Omega)^3$, $M_h\subset{L}_0^2(\Omega)\cap{H}^1(\Omega)$, and $S_h\subset{L_0^2(\Omega)}$ be the conforming Lagrange finite element spaces, and $\bm{D}_h\subset\bm{H}_0(\mathrm{div};\Omega)$ be the $\bm{H}(\mathrm{div})$-conforming finite element space. It is assumed that the pair $\bm{D}_h\times{S}_h$ satisfies the corresponding LBB condition \cite{2021_Ern}: there exists a positive constant $C$ which is independent of $h$, such that
\begin{equation*}
    \inf_{\zeta_h\in{S}_h}\sup_{\bm{K}_h\in\bm{D}_h}\frac{(\nabla\cdot\bm{K}_h,\zeta_h)}{\|\bm{K}_h\|_{\bm{H}(\mathrm{div};\Omega)}\|\zeta_h\|_{L^2(\Omega)}}\ge{C}.
\end{equation*}
We use a uniform partition of the time interval $[0,T]$ with the time step size $\tau$. Then we define
\begin{equation*}
    \delta_\tau\chi^{n+1}=\frac{\gamma_0\chi^{n+1}-\widehat{\chi}}{\tau},
\end{equation*}
where $\chi$ is a generic variable, and
\begin{equation*}
    \widehat{\chi}=
    \begin{cases}
        \chi^n, & k=1, \\
        2\chi^n-\frac{1}{2}\chi^{n-1}, & k=2,
    \end{cases}
    ; \quad \gamma_0=
    \begin{cases}
        1, & k=1, \\
        3/2, & k=2.
    \end{cases}
\end{equation*}
We also recall the extrapolation formulas
\begin{equation*}
    \widetilde{\chi}^{n+1}=
    \begin{cases}
        \chi^n, & k=1, \\
        2\chi^n-\chi^{n-1}, & k=2.
    \end{cases}
\end{equation*}
With these notations, the discrete BDF finite element scheme is given as follows:

\textit{Step} 1: find $\left(\phi_h^{n+1},\mu_h^{n+1}\right)\in{Q}_h\times{Q}_h$, such that for all $(\psi_h,\omega_h)\in{Q}_h\times{Q}_h$,
\begin{equation}
\label{step_1}
    \begin{aligned}
        \big(\delta_\tau\phi_h^{n+1},\omega_h\big)+\frac{1}{\mathrm{Pe}}\big(\nabla\mu_h^{n+1},\nabla\omega_h\big)= & \big(\widetilde{\phi}_h^{n+1}\widetilde{\bm{u}}_h^{n+1},\nabla\omega_h\big) , \\
        \big(\mu_h^{n+1},\psi_h\big)-\mathrm{Cn}^2\big(\nabla\phi_h^{n+1},\nabla\psi_h\big)-\left(\phi_h^{n+1},\psi_h\right)=  & \big((\widetilde{\phi}_h^{n+1})^3-2\widetilde{\phi}_h^{n+1},\psi_h\big).
    \end{aligned}
\end{equation}
Then we evaluate $\rho_h^{n+1}$, $\eta_h^{n+1}$ and $\sigma_h^{n+1}$ by (exemplarily for $\rho_h^{n+1}$)
\begin{equation*}
    \rho_h^{n+1}=\frac{\rho_+-\rho_-}{2\rho_r}\mathscr{C}\left(\phi_h^{n+1}\right)+\frac{\rho_++\rho_-}{2\rho_r},
\end{equation*}
where we define
\begin{equation*}
    \mathscr{C}\left(\phi_h^{n+1}\right)=
    \begin{cases}
        \phi_h^{n+1}, & \text{if}\;\left|\phi_h^{n+1}\right|\le{1}, \\
        \mathrm{sign}(\phi_h^{n+1}), & \text{otherwise}.
    \end{cases}
\end{equation*}

\textit{Step} 2: find $\bm{u}_h^{n+1}\in\bm{X}_h$, such that for all $\bm{v}_h\in\bm{X}_h$,
\begin{equation}
\label{step_2}
    \begin{aligned}
        \big(\rho_h^{n+1}\delta_\tau\bm{u}_h^{n+1},\bm{v}_h\big)+\frac{2}{\mathrm{Re}}\big(\eta_h^{n+1}\D(\bm{u}_h^{n+1}),\D(\bm{v}_h)\big)= & -\big((\rho_h^{n+1}\widetilde{\bm{u}}_h^{n+1}-\rho_d\frac{1}{\mathrm{Pe}}\nabla\mu_h^{n+1})\cdot\nabla\widetilde{\bm{u}}_h^{n+1},\bm{v}_h\big) \\
        & +\big(\widetilde{p}_h^{n+1},\nabla\cdot\bm{v}_h\big)+\frac{1}{\mathrm{We}\mathrm{Cn}}\big(\mu_h^{n+1}\nabla\phi_h^{n+1},\bm{v}_h\big)\\
        & +\mathrm{N}\big(\widetilde{\bm{J}}_h^{n+1}\times\bm{B},\bm{v}_h\big)-\frac{1}{\mathrm{Fr}}\big(\rho_h^{n+1}\bm{e}_3,\bm{v}_h\big).
    \end{aligned}
\end{equation}

\textit{Step} 3: find $p_h^{n+1}\in{M}_h$, such that for all $q_h\in{M}_h$,
\begin{equation}
\label{step_3}
    \big(\nabla({p}_h^{n+1}-{p}_h^n),\nabla{q}_h\big)=-\frac{\gamma_0\vartheta}{\tau}\big(\nabla\cdot\bm{u}_h^{n+1},q_h\big),
\end{equation}
with $\vartheta=\min\left(\rho_+,\rho_-\right)/\rho_r$.

\textit{Step} 4: find $V_h^{n+1}\in{M}_h$, such that for all $\Lambda_h\in{M}_h$
\begin{equation}
\label{step_4}
    \left(\sigma_h^{n+1}\nabla{V}_h^{n+1},\nabla\Lambda_h\right)=\left(\sigma_h^{n+1}\bm{u}_h^{n+1}\times\bm{B},\nabla\Lambda_h\right).
\end{equation}

\textit{Step} 5: find $(\bm{J}_h^{n+1},\xi_h^{n+1})\in\bm{D}_h\times{S}_h$, such that for all $(\bm{K}_h,\zeta_h)\in\bm{D}_h\times{S}_h$
\begin{equation}
\label{step_5}
    \begin{aligned}
        \left(\bm{J}_h^{n+1},\bm{K}_h\right)-\left(\xi_h^{n+1},\nabla\cdot\bm{K}_h\right) & =\left(\sigma_h^{n+1}(-\nabla{V}_h^{n+1}+\bm{u}_h^{n+1}\times\bm{B}),\bm{K}_h\right), \\
        -\left(\nabla\cdot\bm{J}_h^{n+1},\zeta_h\right) & =0,
    \end{aligned}
\end{equation}
where $\xi_h^{n+1}$ is a Lagrange multiplier.

In the above scheme, we have exploited a linearizing technique in \cite{2015_Shen} and the pressure stabilization method developed by Guermond and Salgado \cite{2009_Guermond}. Additionally, we solve the electrostatic potential $V$ through the electric potential Poisson equation \cite{2014_Zhang} and the current density $\bm{J}$ is then projected into the $\bm{H}(\mathrm{div})$-conforming finite element spaces. These techniques greatly facilitate numerical computations. For time-independent coefficient matrix in equations \eqref{step_1} and \eqref{step_5}, expensive block type preconditioners can be employed, and for Poisson-type equations in equations \eqref{step_2}-\eqref{step_4}, robust multilevel preconditioners can be utilized. While this scheme offers a straightforward implementation, it does not hold a discrete energy law, which is particularly desirable in diffuse interface computations \cite{2015_Shen, 2024_Wang}. Despite this restriction, one can readily verify that the scheme \eqref{step_1}-\eqref{step_5} admits a unique solution, and is globally mass conservative and charge conservative (we refer to \cite{2024_Wang} for some details),
\begin{equation*}
    \int_\Omega\phi_h^{n+1}\mathrm{d}\bm{x}=\int_\Omega\phi_h^n\mathrm{d}\bm{x}
    \qquad \mbox{and} \qquad
    \nabla\cdot\bm{J}_h=0.
\end{equation*}

\begin{remark}
     Via its definition, the order parameter $\phi$ should physically belong to $[-1,1]$, yet the Ginzburg-Landau energy cannot guarantee this bound at a theoretical level, unlike the Flory-Huggins singular (logarithmic) energy \cite{2019_Giorgini}. In computations, we employ a cut-off approach to numerically preserve this physical constraint, which is a common strategy in diffuse interface simulations using Ginzburg-Landau energy \cite{2023_Khanwale, 2015_Shen, 2024_Wang}. This technique enables stable simulations with large parameter ratios.
\end{remark}
\begin{remark}
     Practical MHD engineering applications usually involves large parameter ratios, with a typical order $\mathcal{O}(10^3)$ for density, $\mathcal{O}(10^2)$ for dynamic viscosity, and $\mathcal{O}(10^7)$ for electrical conductivity; see \cite{2024_Gou, 2014_Zhang} and the references therein for more details in laboratory experiments. This results in extremely large Reynolds numbers and Stuart numbers, in turn demanding a highly robust solver to resolve the flow dynamics. In addition, many MHD processes exhibit multiscale features \cite{2019_Herreman, 2018_Thomas, 2018_Tucs}, requiring extremely small interface thicknesses to achieve a better simulation result. Undoubtedly, this will further consume a vast amount of computational resources. Under such a circumstance, developing an asymptotic-preserving and robust solver for realistic MHD simulations is essential and is of our future interest.
\end{remark}

\subsection{The magnetic damping effects on bubble dynamics}

\begin{figure}[H]
	\centering
	\subfloat[Horizontal magnetic field]{
		\includegraphics[width=0.25\linewidth]{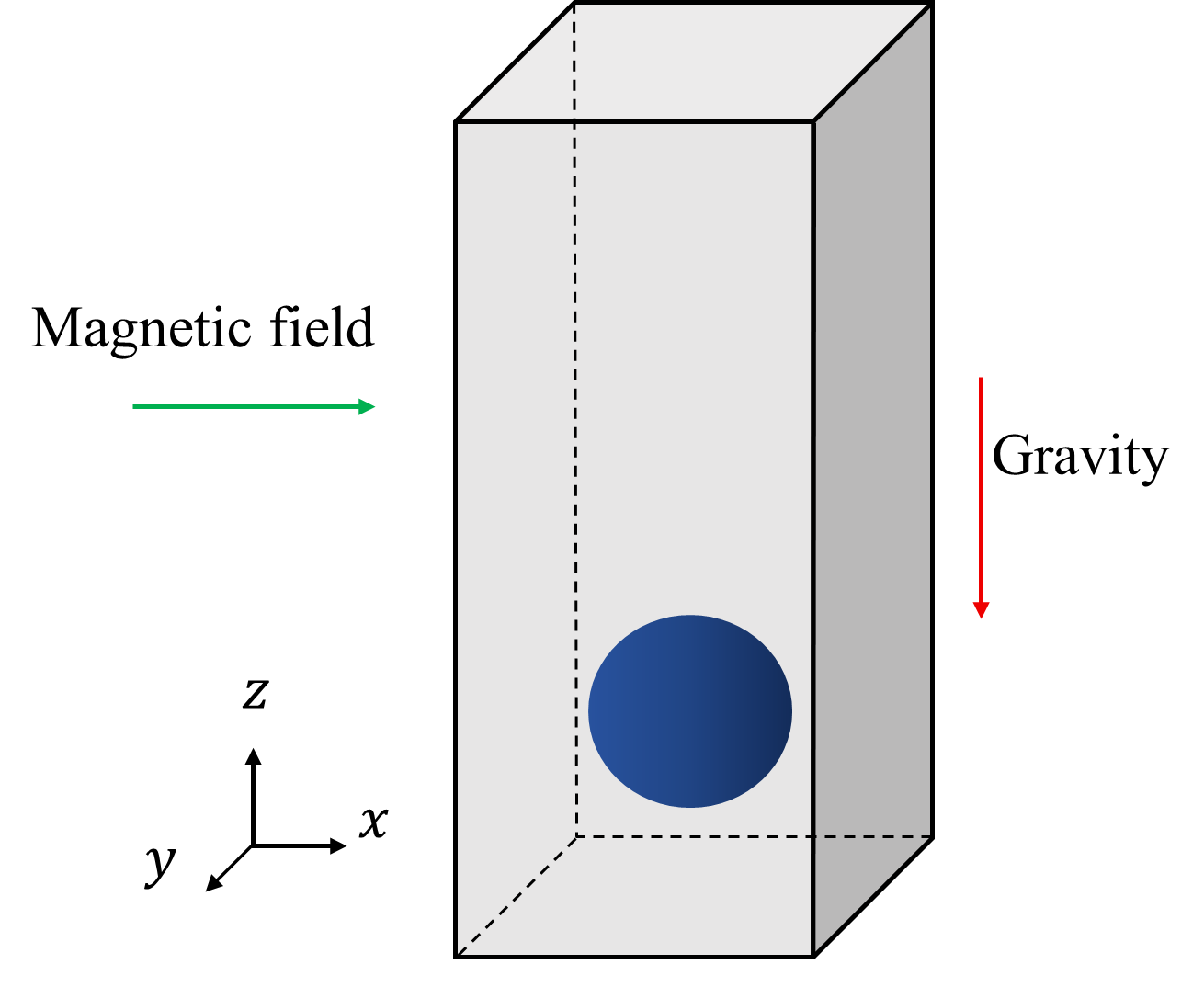}
	}
	\hspace{0.1\linewidth}
	\subfloat[Vertical magnetic field]{
		\includegraphics[width=0.25\linewidth]{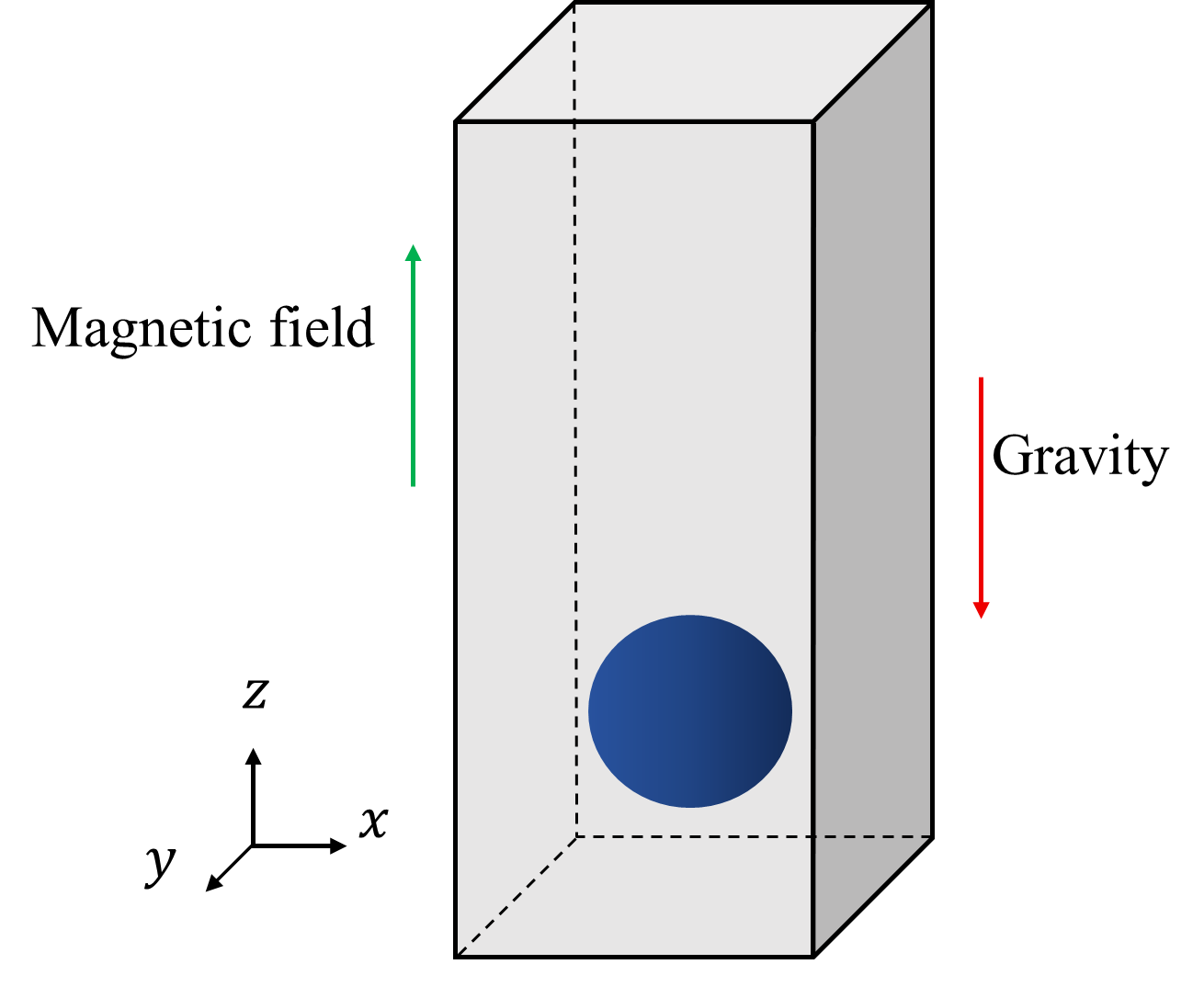}
	}
	\caption{Setup of single rising bubble under magnetic fields.}
	\label{problem_set_up}
\end{figure}

In this section, we present several three-dimensional numerical examples to demonstrate the magnetic damping effects on bubble dynamics. An illustration can be seen in Figure \ref{problem_set_up}. This problem finds its applications in metallurgy processes \cite{2018_Thomas}, where injected bubbles serve to stir and homogenize molten metals, with magnetic fields providing a non-intrusive approach to control the motion of bubbles. Meanwhile, it has been extensively investigated both experimentally \cite{2024_Gou} and numerically \cite{2014_Zhang, 2014_Zhang_pof, 2016_Zhang}. Here, we focus on the scenario of a single bubble rising in a liquid column under either an external horizontal or vertical magnetic field. In this benchmark problem, a bubble of a lighter fluid with radius $0.5\;\mathrm{m}$ is centered at $(0.5\;\mathrm{m},0.5\;\mathrm{m},0.5\;\mathrm{m})$ in the domain $\Omega=(0,1\;\mathrm{m})\times(0,1\;\mathrm{m})\times(0,2\;\mathrm{m})$, corresponding to the following initial profile for $\phi$,
\begin{equation*}
    \phi^0=\tanh\left(\frac{\sqrt{(x-0.5)^2+(y-0.5)^2+(z-0.5)^2}-0.25}{\sqrt{2}\epsilon}\right).
\end{equation*}
We will use the international system of units and ignore the units later if there is no confusion. The relevant physical parameters are
\begin{itemize}
    \item densities ($\mathrm{kg}\cdot\mathrm{m}^{-3}$): $\rho_+=1000$, $\rho_-=1$;
    \item dynamic viscosities ($\mathrm{kg}\cdot\mathrm{m}^{-1}\cdot\mathrm{s}^{-1}$): $\eta_+=10$, $\eta_-=0.1$;
    \item electrical conductivities ($\mathrm{S}\cdot\mathrm{m}^{-1}$): $\sigma_+=1000$, $\sigma_-=1$;
    \item surface tension ($\mathrm{N}\cdot\mathrm{m}^{-1}$): $\lambda=1.96$;
    \item gravitational acceleration ($\mathrm{m}\cdot\mathrm{s}^{-2}$): $g=0.98$;
    \item magnetic fields ($\mathrm{T}$): a horizontal magnetic field $\bm{B}=(B_r,0,0)^\top$ and a vertical magnetic field $\bm{B}=(0,0,B_r)^\top$, with $B_r=3$, $5$, and $7$.
\end{itemize}

Our numerical implementation relies on the \texttt{deal.II} finite element library \cite{2023_Arndt} with \texttt{ParaView} post-processing software (\url{www.paraview.org}). The proposed BDF2 scheme, initialized via BDF1, is utilized in simulation. Additionally, we discretize $Q_h$, $\bm{X}_h$, and $M_h$ with the continuous $\mathbb{Q}_1$ element, $S_h$ with the discontinuous $\mathbb{DGQ}_0$ element, and $\bm{D}_h$ with the lowest order Raviart-Thomas $\mathbb{RT}_0$ element. To solve the resulting large-scale algebraic systems, we employ the flexible generalized minimum residual (FGMRES) algorithm with preconditioning: additive Schwarz methods for equations \eqref{step_1} and \eqref{step_5}, and algebraic multigrid (AMG) methods for equations \eqref{step_2}-\eqref{step_4}. Additionally, we use a iteration tolerance of $\varepsilon_{iter}=10^{-10}$ with the initial guess set to the extrapolated value.

\begin{figure}[H]
    \centering
    \includegraphics[width=0.3\linewidth]{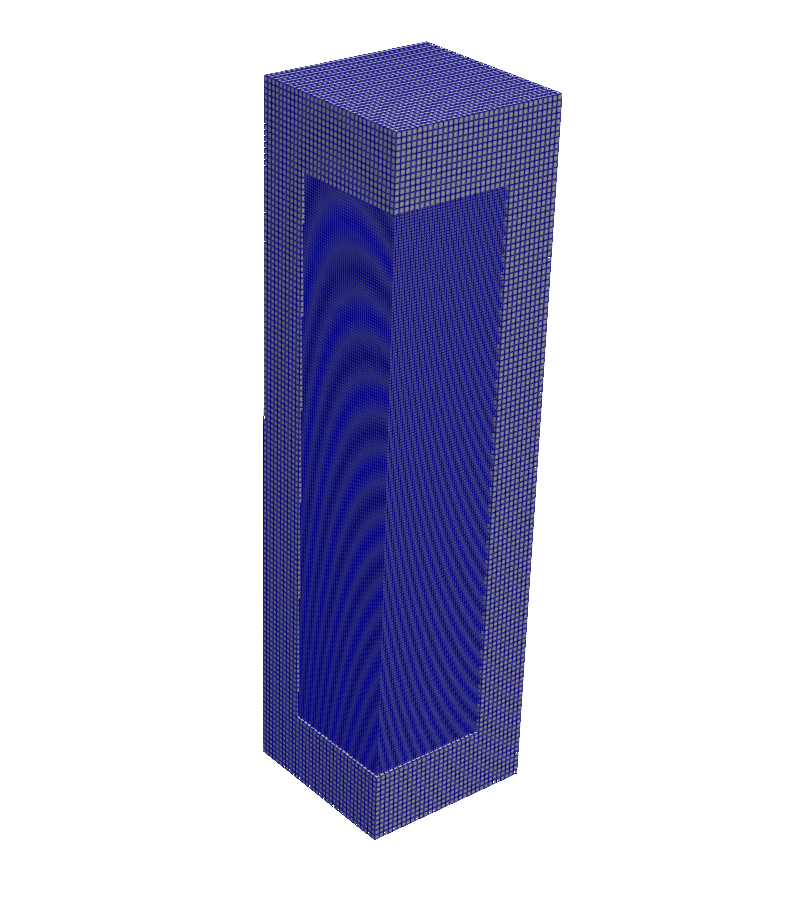}
    \caption{Mesh, quarter of domain.}
    \label{mesh}
\end{figure}

While adaptive mesh refinement (AMR) techniques lead to a poor mass conservation result \cite{2023_Khanwale} and present challenges in designing a desirable refinement strategy, we only adopt a fixed locally refined mesh to resolve the interfacial dynamics. Specifically, we set mesh size $h=2^{-7}$ in subdomain $(0.15,0.85)\times(0.15,0.85)\times(0.2,1.8)$ and $h=2^{-6}$ elsewhere; see Figure \ref{mesh}. This configuration generates 1,906,592 cells, resulting in 3,920,510, 5,880,765, 1,960,255, 1,960,255, and 7,701,408 degrees of freedom for equations \eqref{step_1}-\eqref{step_5}, respectively. In terms of temporal discretization, we set the final dimensional time $T^*=3$ with a fixed time step of $1\times{10}^{-3}$, giving 3,000 total steps. With regard to characteristic quantities, we select $\epsilon=0.005$, $L_r=1$, $u_r=\sqrt{gL_r}$, $\rho_r=\rho_+$, $\eta_r=\eta_+$, and $\sigma_r=\sigma_+$. All results were obtained in a machine with AMD EPYC 9654 and 768 GB RAM.

\begin{figure}[H]
    \centering
    \subfloat[$t^*=1$]{
        \includegraphics[width=0.3\linewidth]{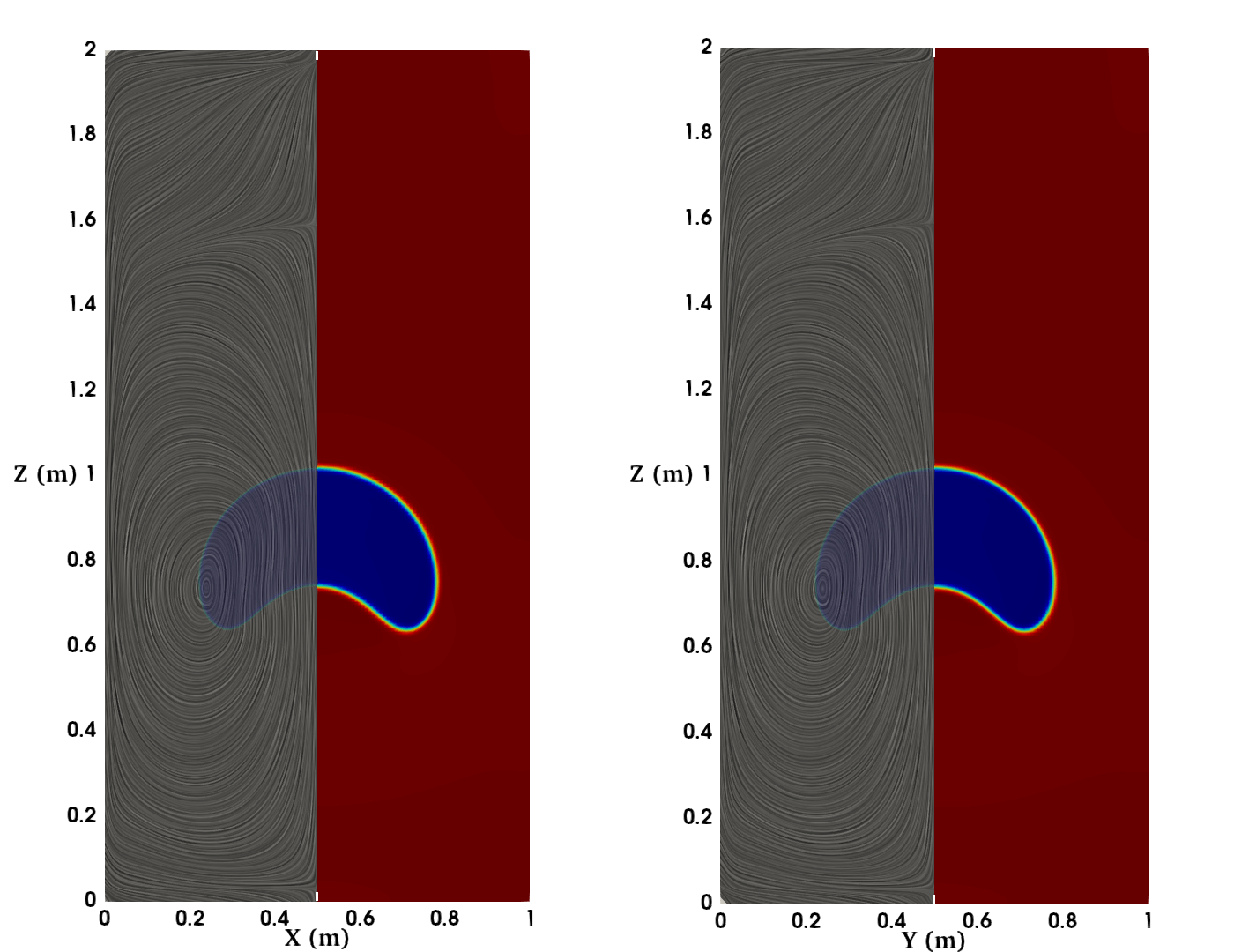}
    }
    \subfloat[$t^*=2$]{
        \includegraphics[width=0.3\linewidth]{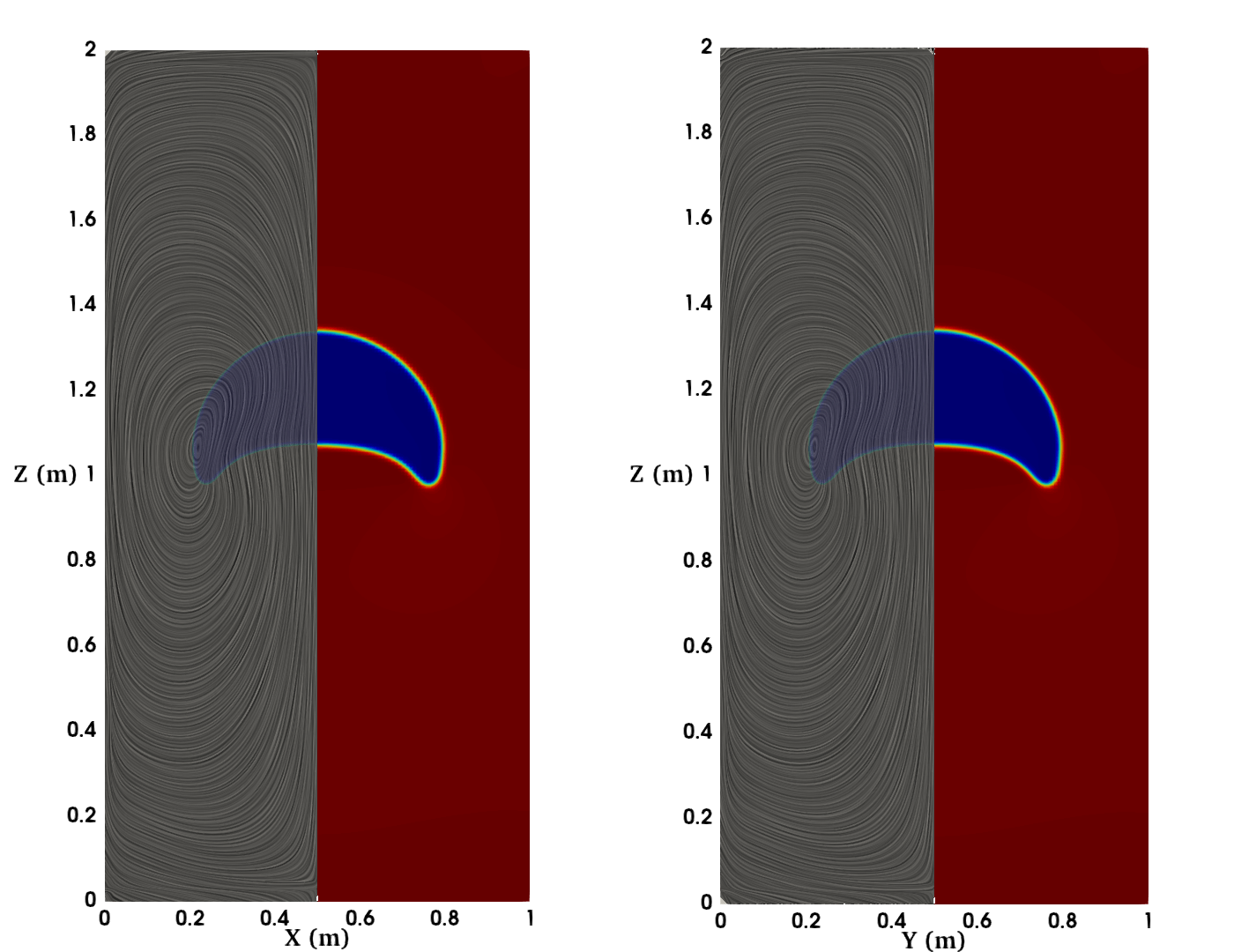}
    }
    \subfloat[$t^*=3$]{
        \includegraphics[width=0.3\linewidth]{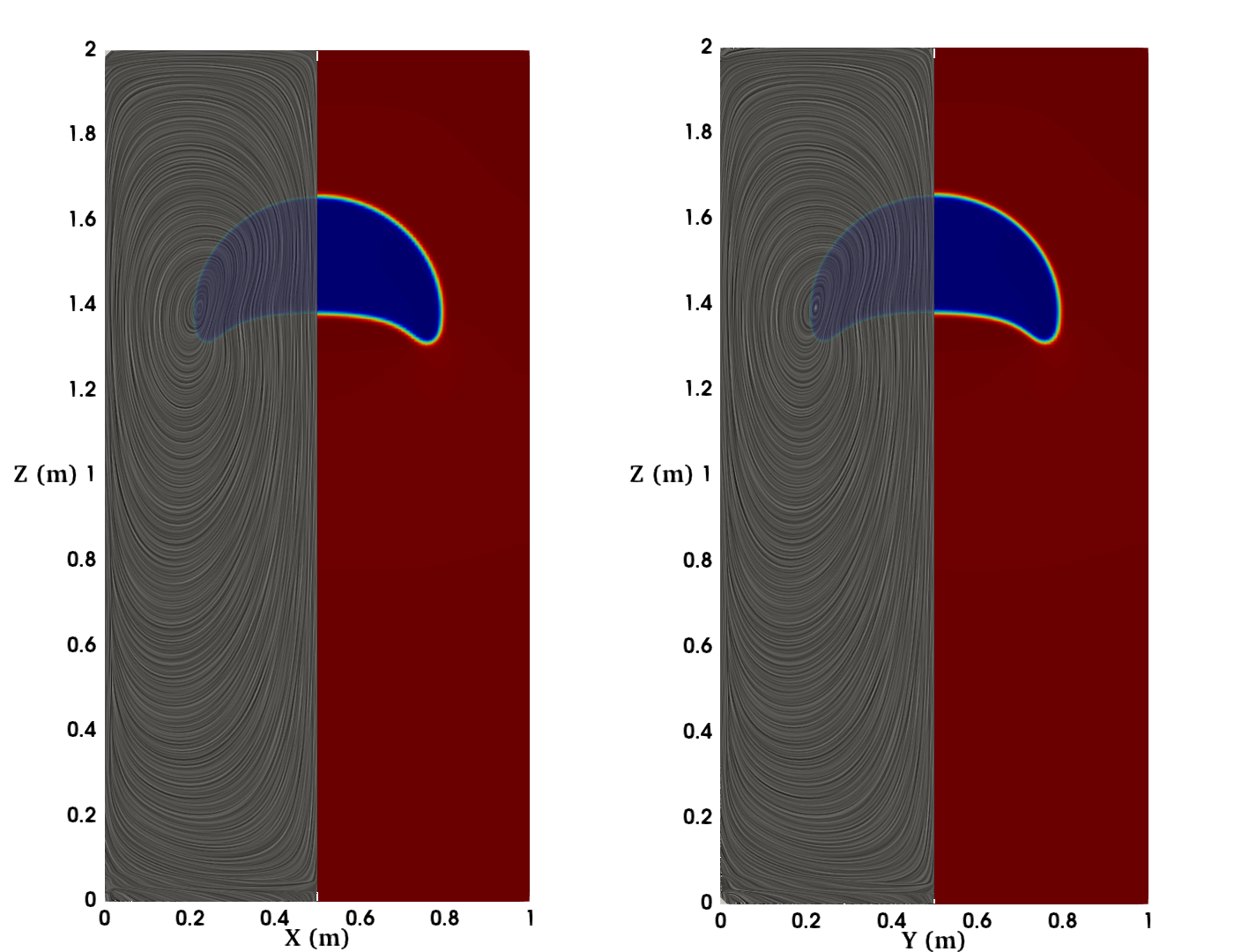}
    }
    \caption{Cross-sections of $\phi_h^n$ without magnetic fields, with streamline distributions of $\bm{u}_h^n$ depicted on the left-hand side.}
    \label{bubble_ref}
\end{figure}

\begin{figure}[H]
    \centering
    \subfloat[$t^*=1$]{
        \includegraphics[width=0.3\linewidth]{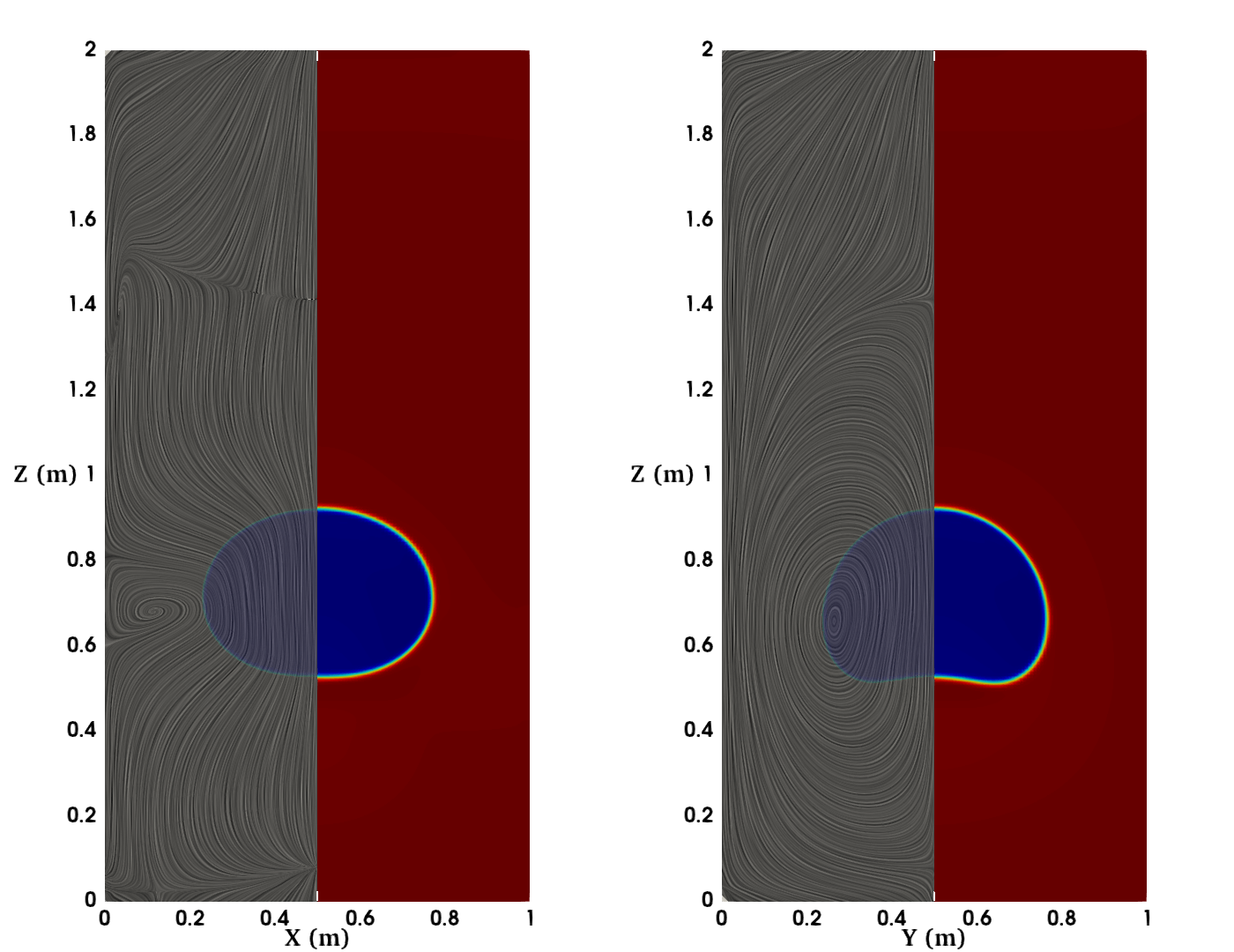}
    }
    \subfloat[$t^*=2$]{
        \includegraphics[width=0.3\linewidth]{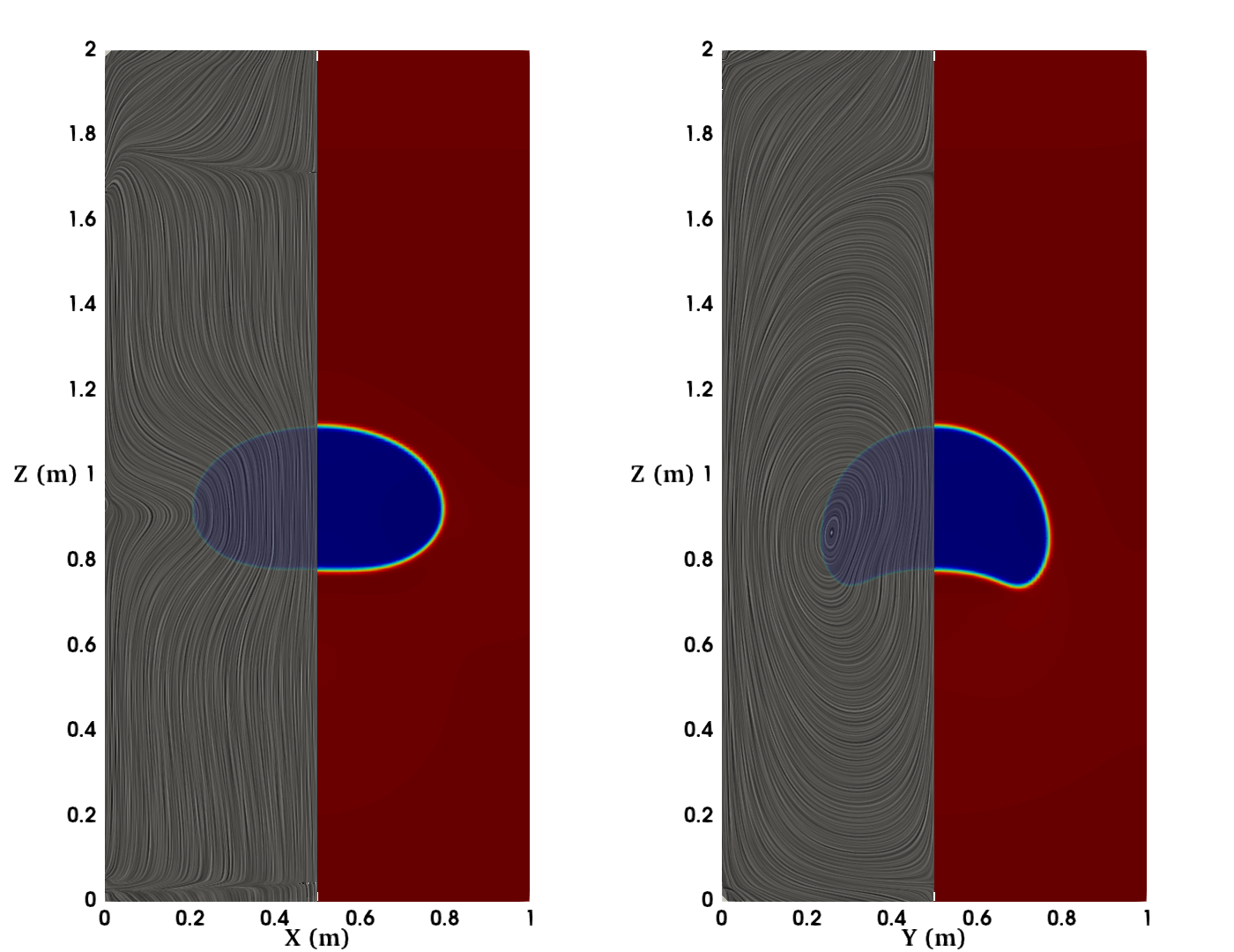}
    }
    \subfloat[$t^*=3$]{
        \includegraphics[width=0.3\linewidth]{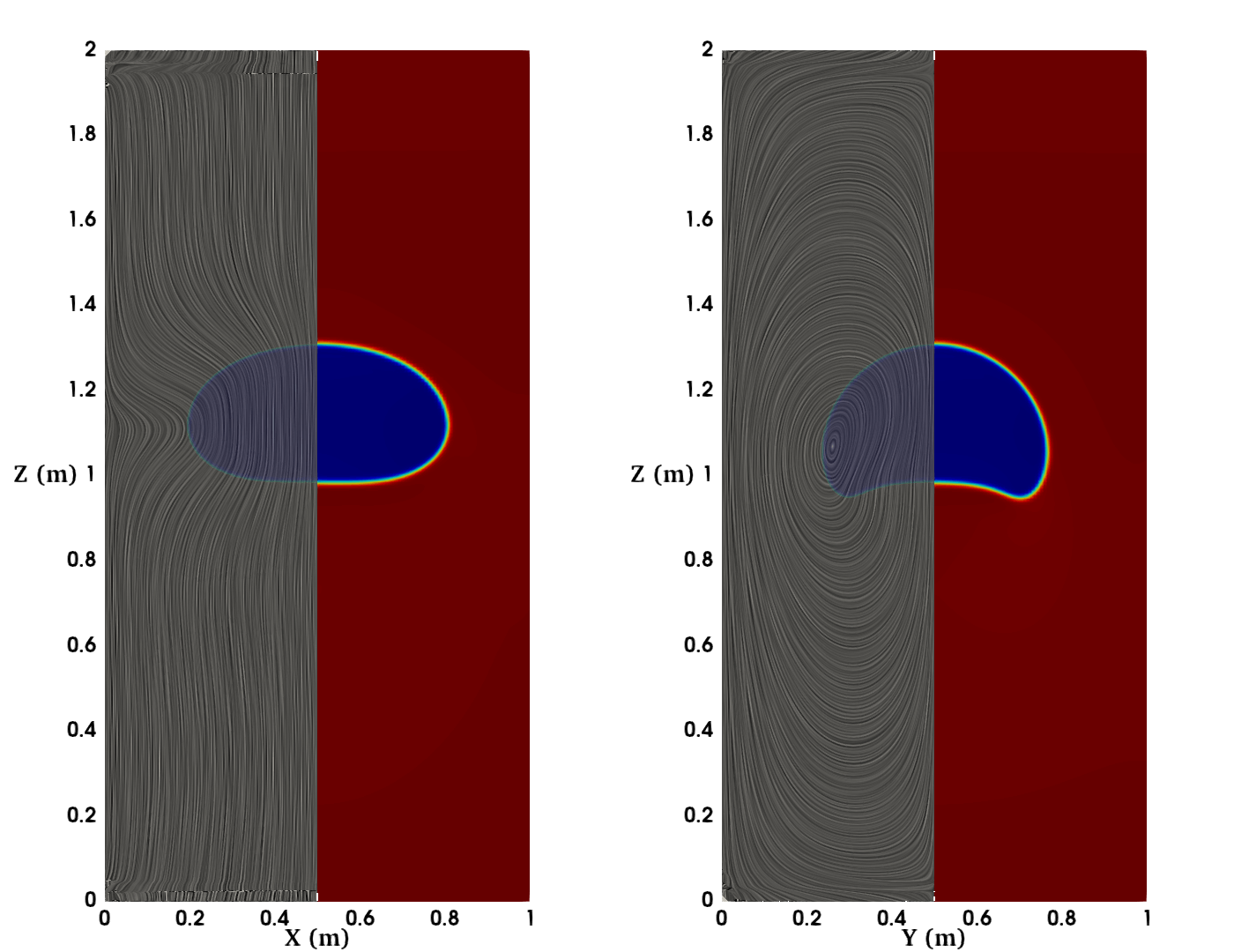}
    }
    \\
    \subfloat[$t^*=1$]{
        \includegraphics[width=0.3\linewidth]{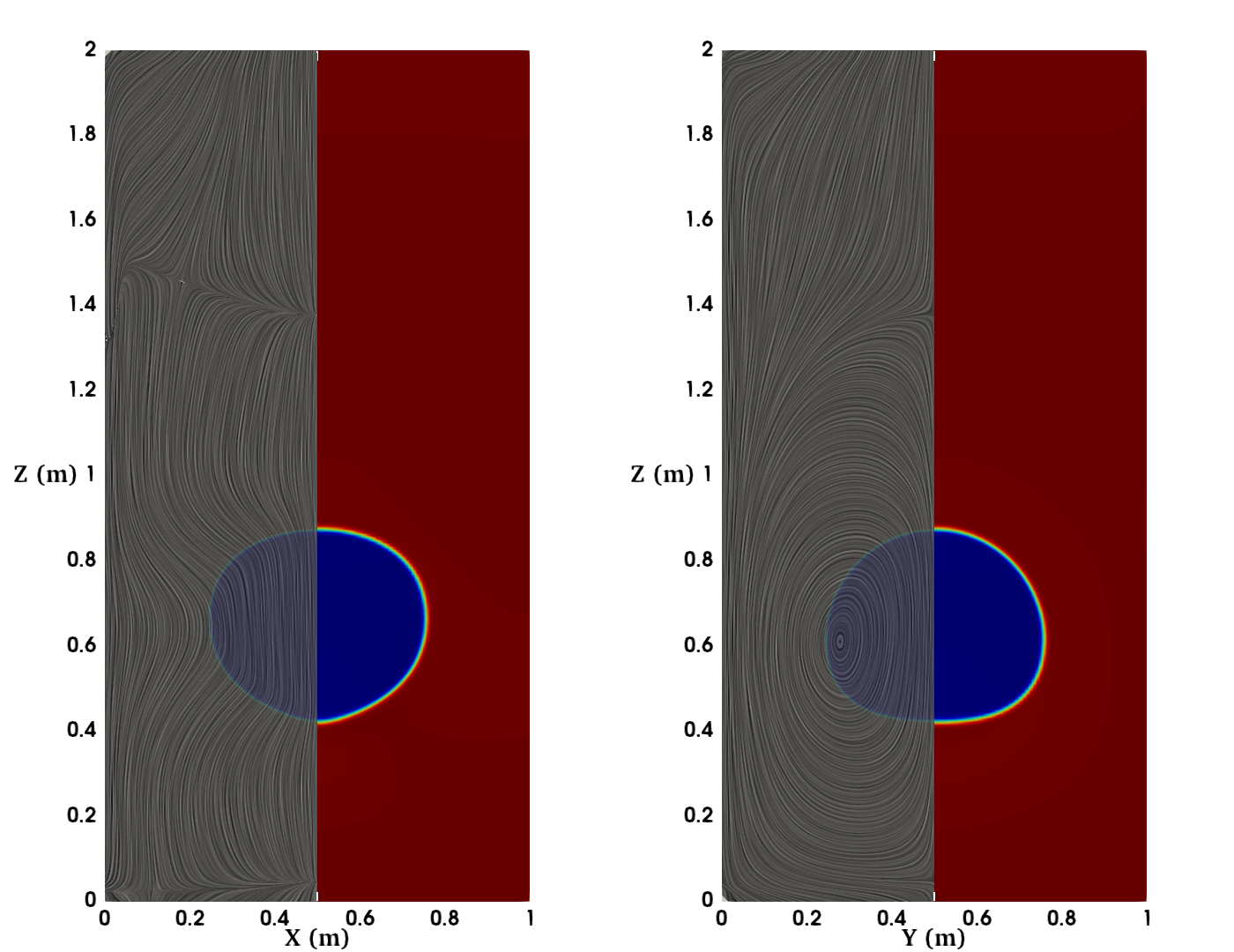}
    }
    \subfloat[$t^*=2$]{
        \includegraphics[width=0.3\linewidth]{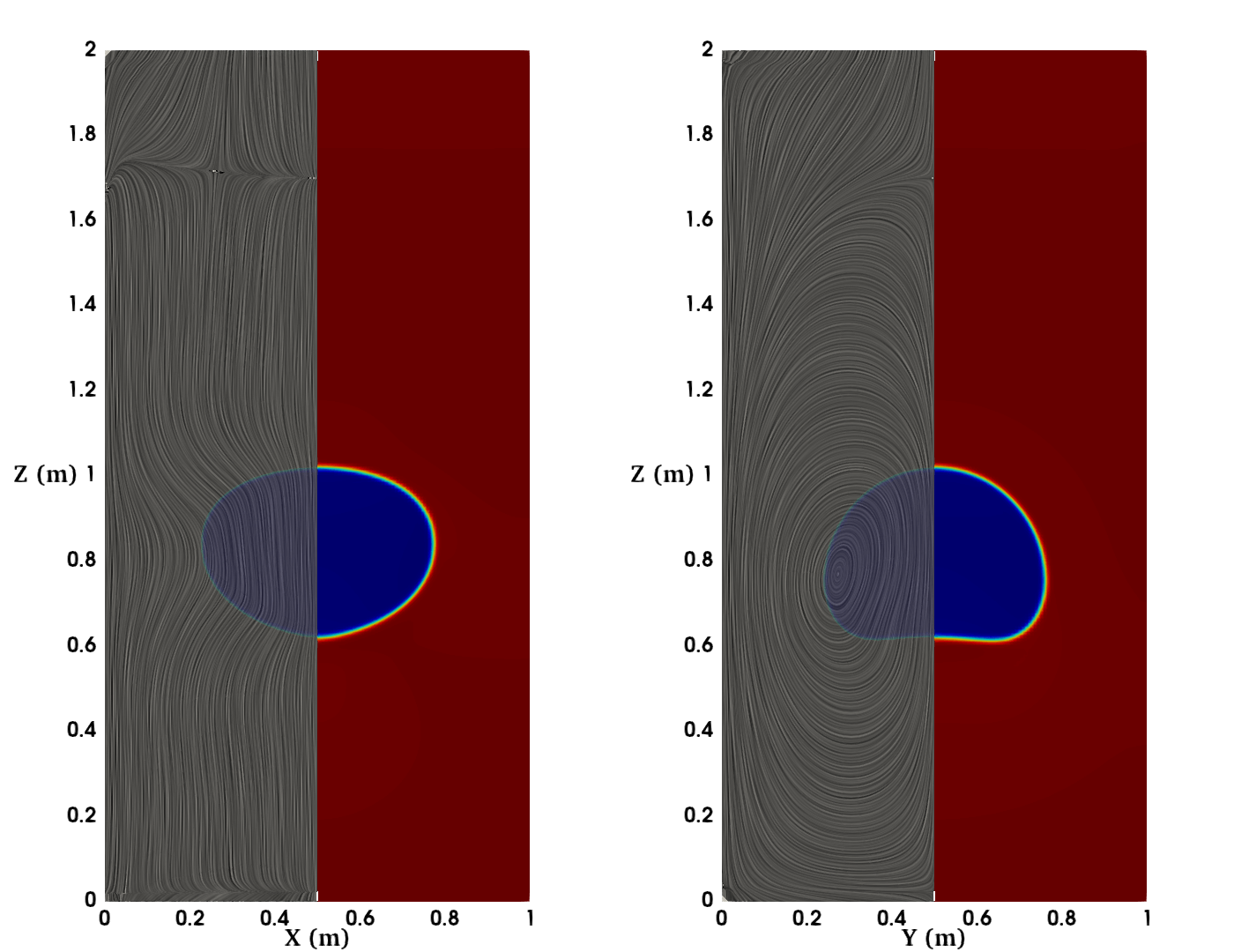}
    }
    \subfloat[$t^*=3$]{
        \includegraphics[width=0.3\linewidth]{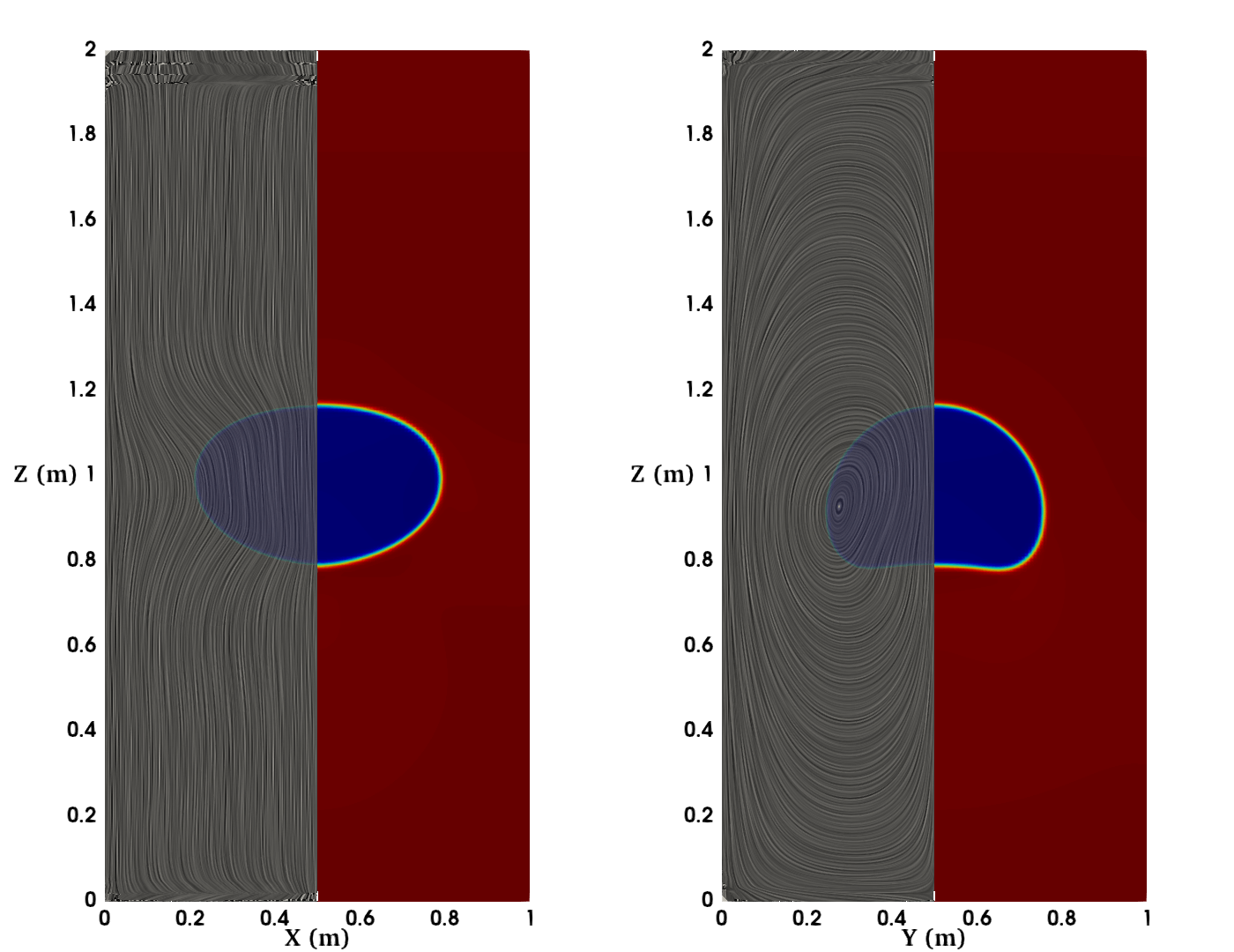}
    }
    \\
    \subfloat[$t^*=1$]{
        \includegraphics[width=0.3\linewidth]{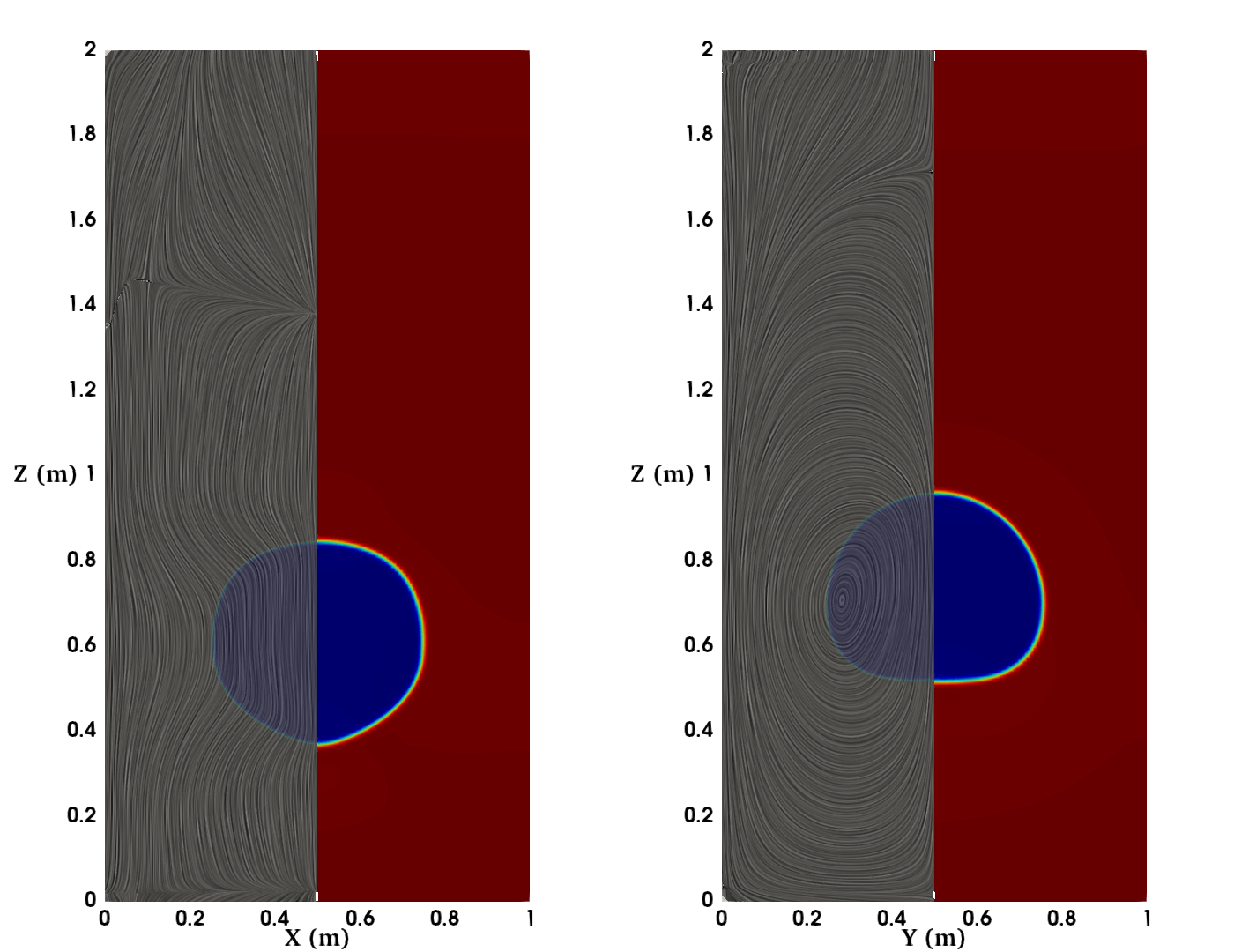}
    }
    \subfloat[$t^*=2$]{
        \includegraphics[width=0.3\linewidth]{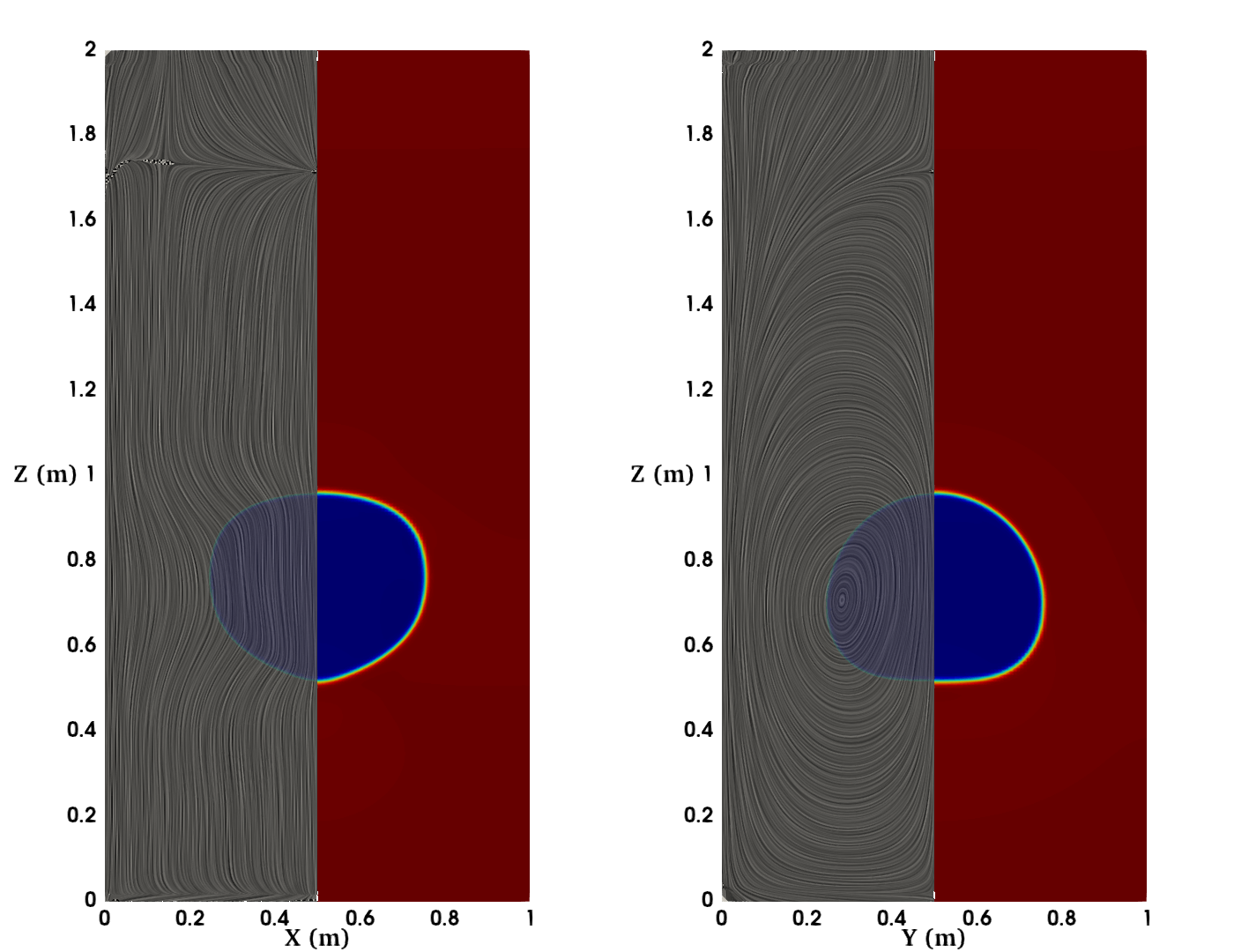}
    }
    \subfloat[$t^*=3$]{
        \includegraphics[width=0.3\linewidth]{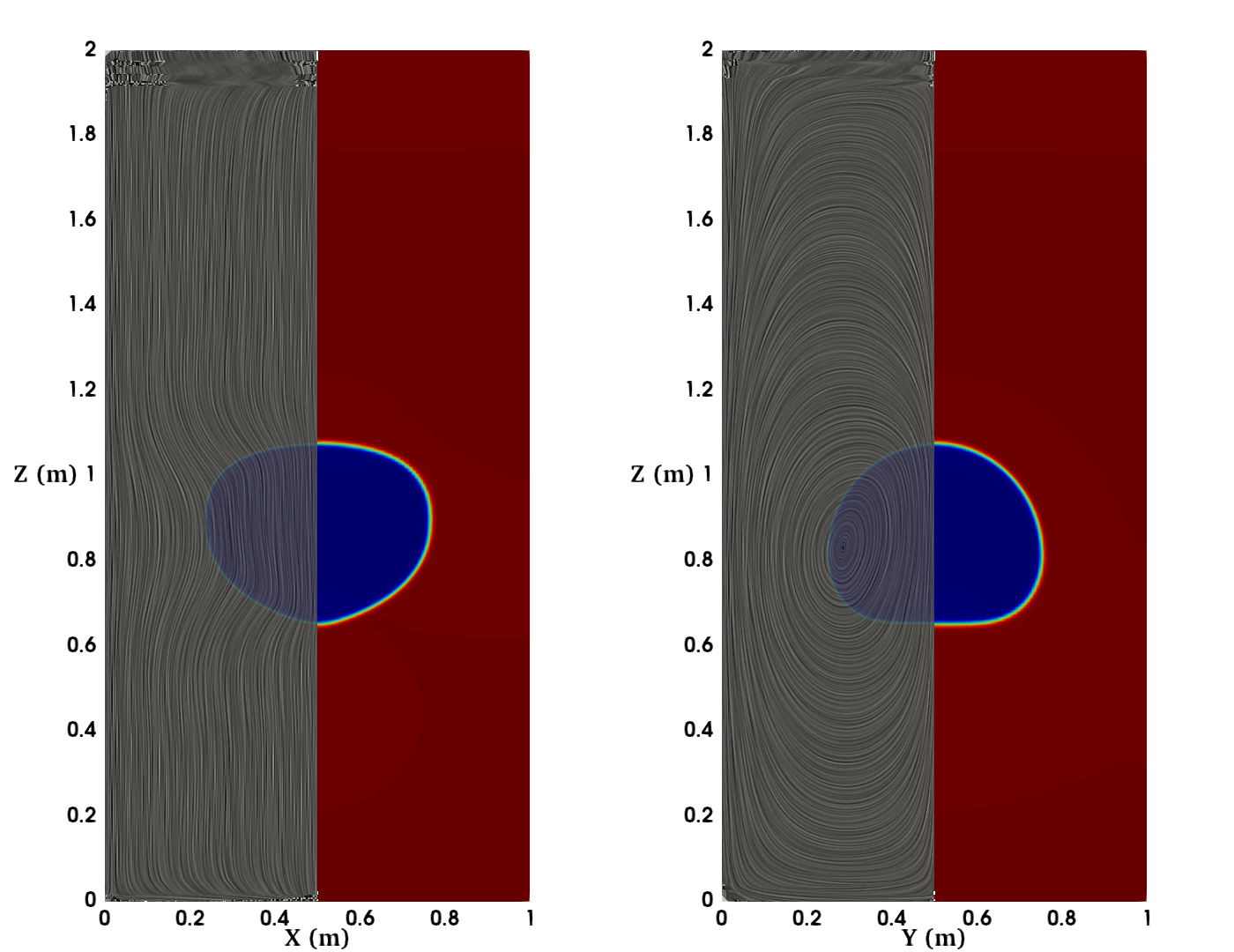}
    }
    \caption{Cross-sections of $\phi_h^n$ under horizontal magnetic fields, with streamline distributions of $\bm{u}_h^n$ depicted on the left-hand side. Top row: $B_r=3$; middle row: $B_r=5$; bottom row: $B_r=7$.}
    \label{bubble_hor}
\end{figure}

\begin{figure}[H]
    \centering
    \subfloat[$t^*=1$]{
        \includegraphics[width=0.3\linewidth]{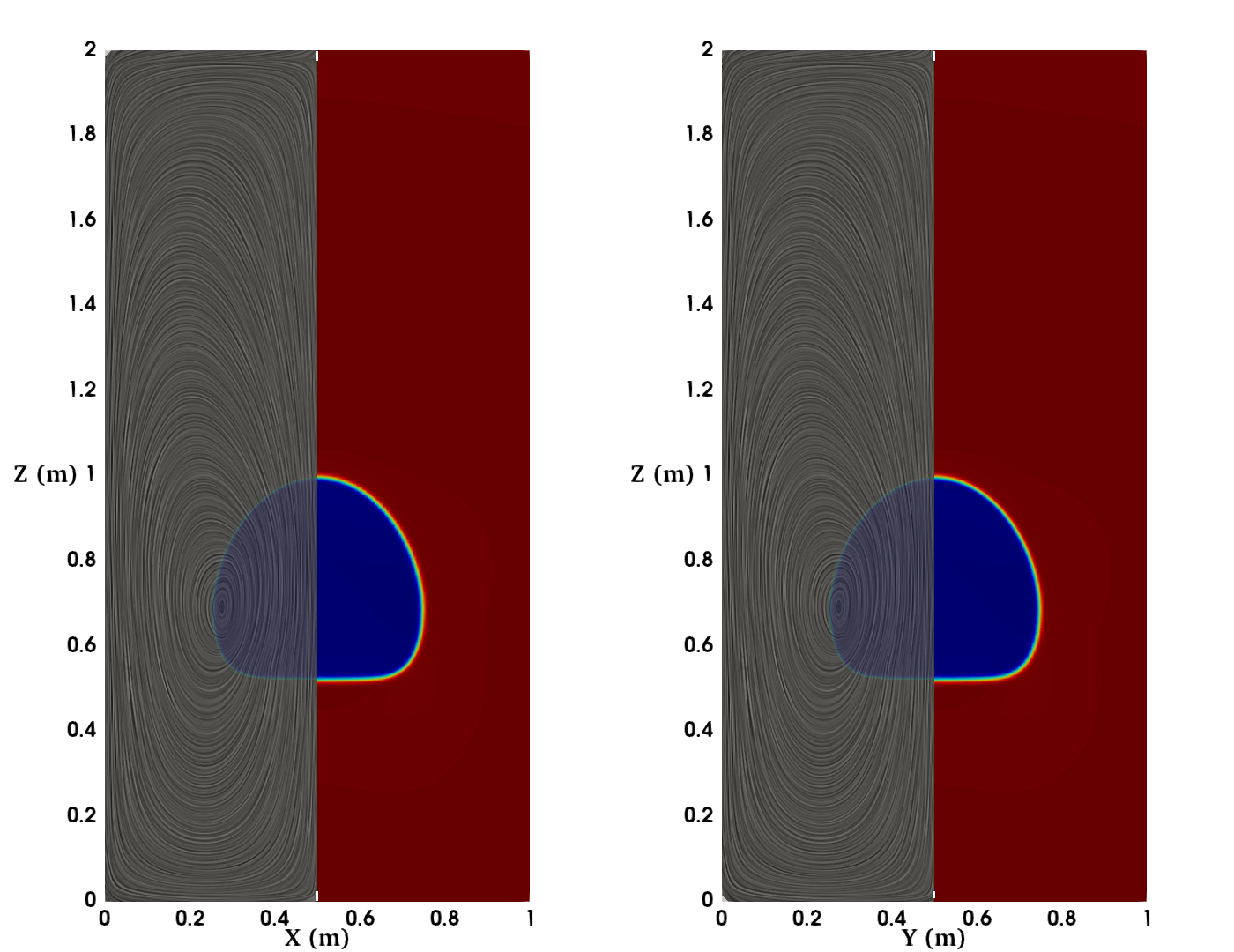}
    }
    \subfloat[$t^*=2$]{
        \includegraphics[width=0.3\linewidth]{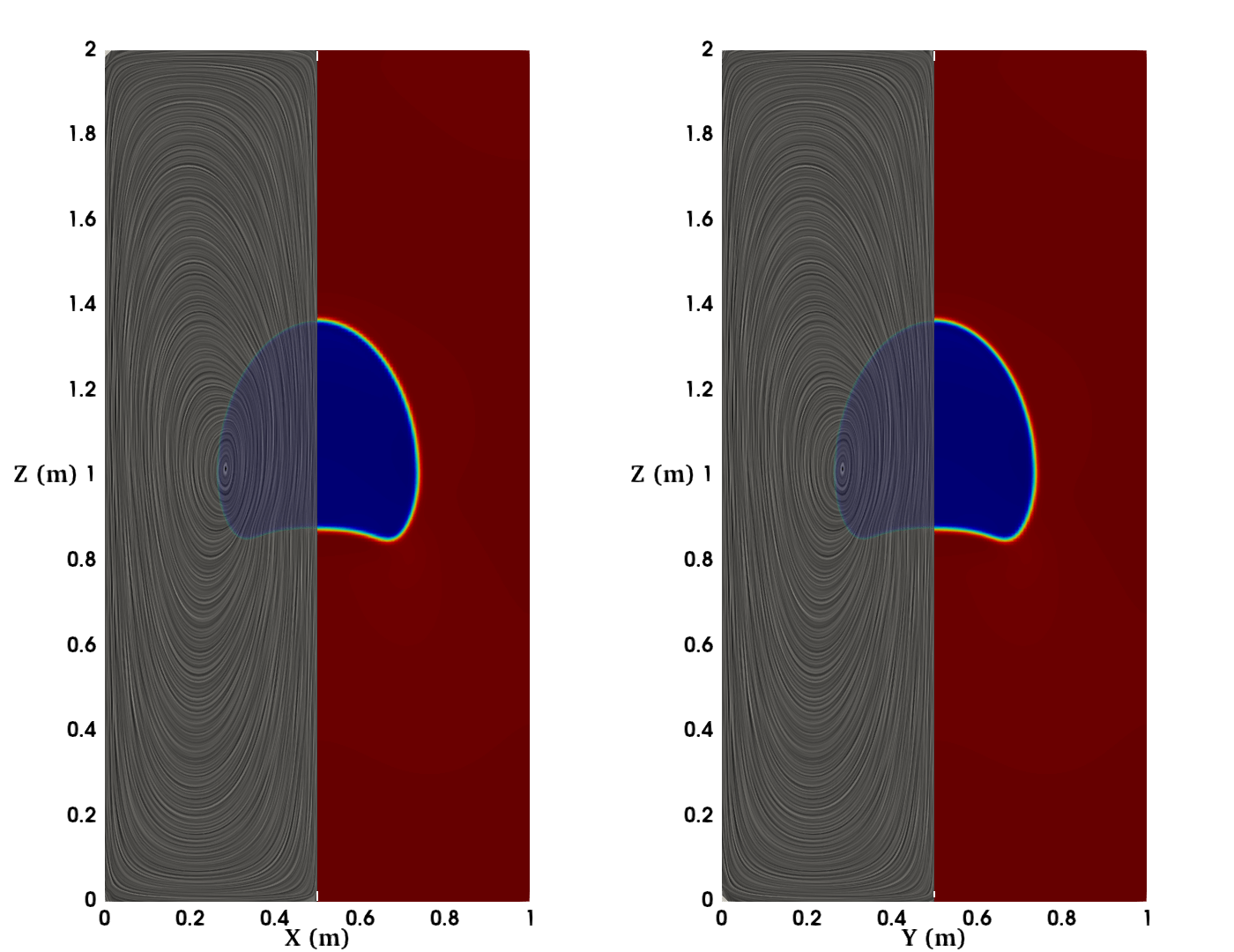}
    }
    \subfloat[$t^*=3$]{
        \includegraphics[width=0.3\linewidth]{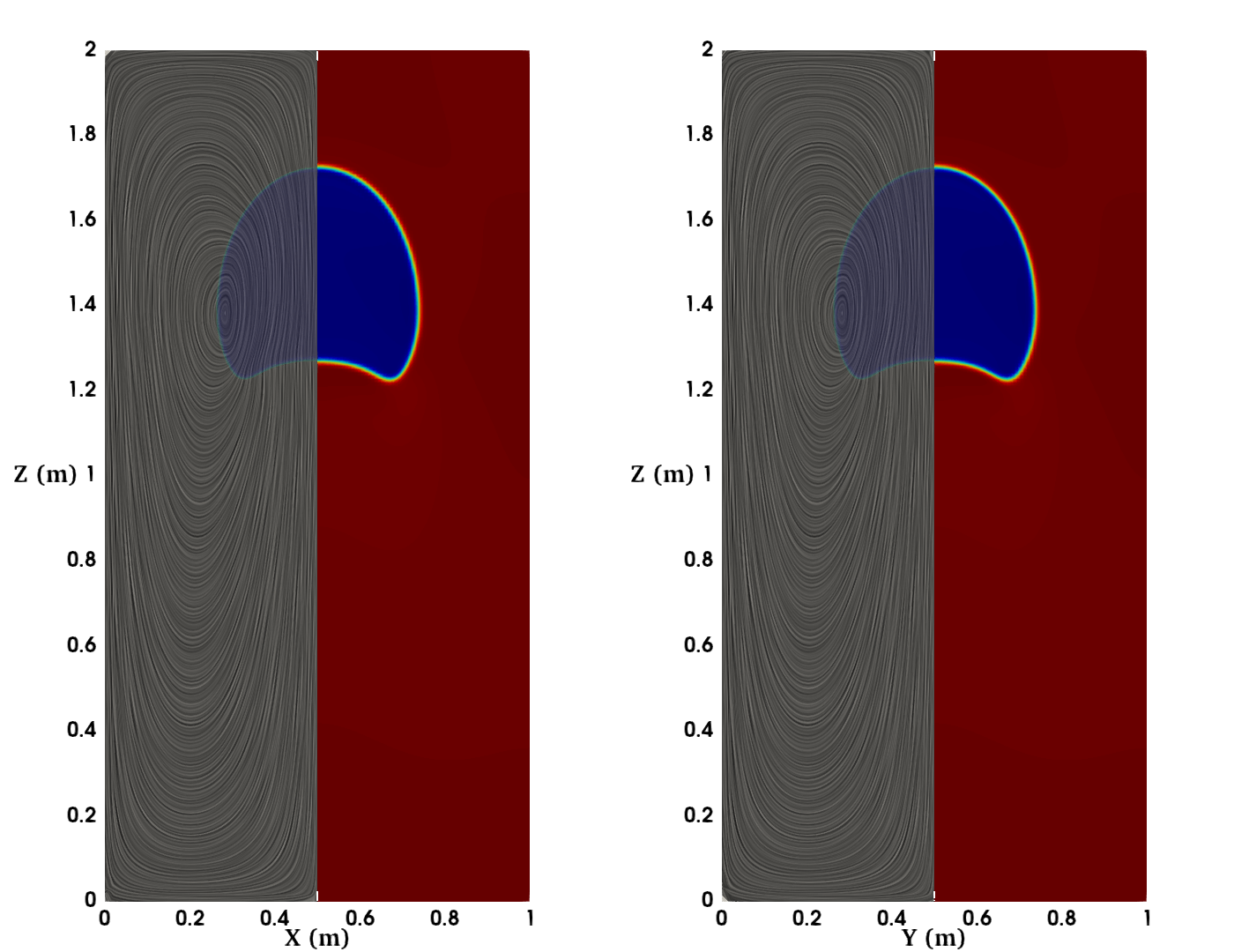}
    }
    \\
    \subfloat[$t^*=1$]{
        \includegraphics[width=0.3\linewidth]{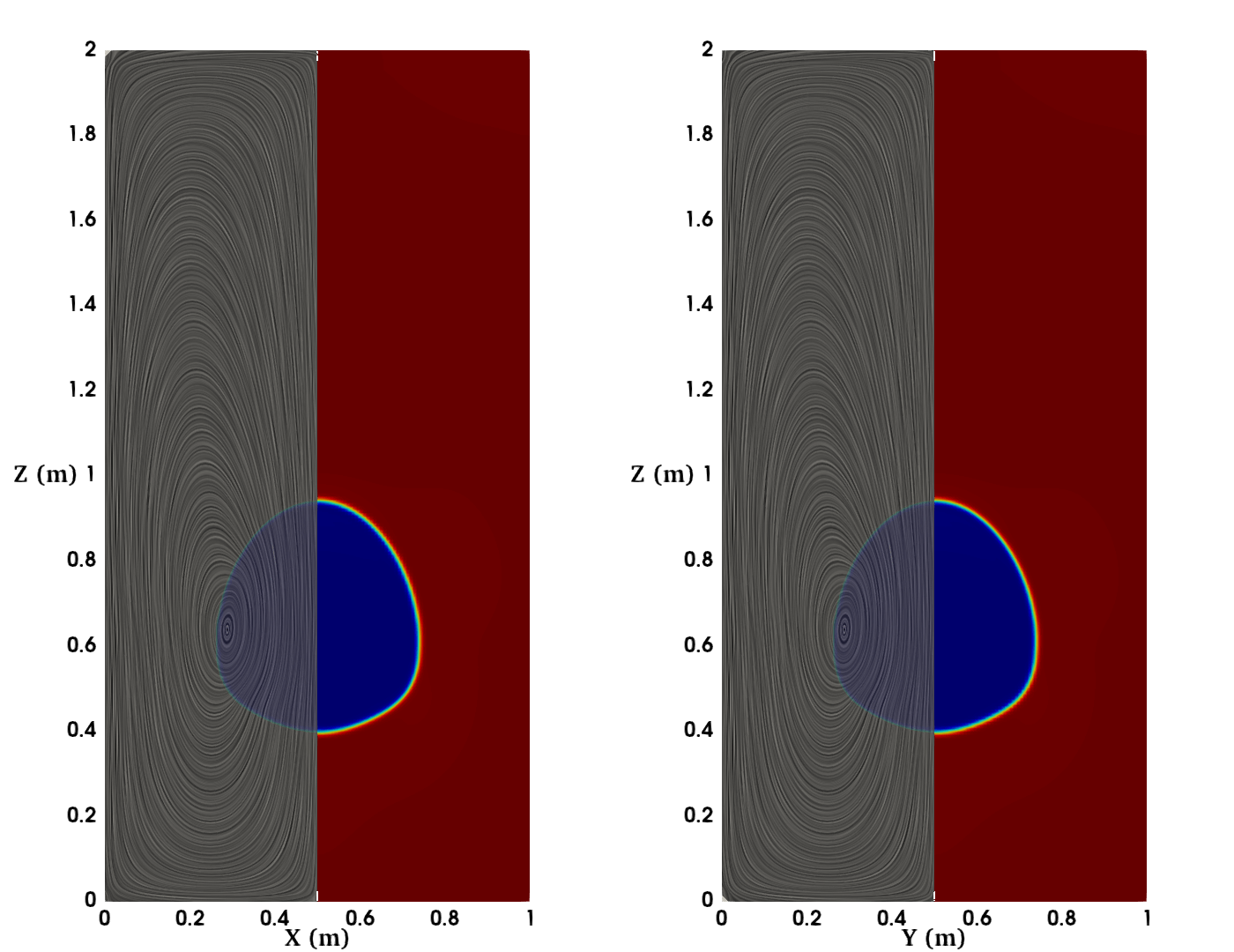}
    }
    \subfloat[$t^*=2$]{
        \includegraphics[width=0.3\linewidth]{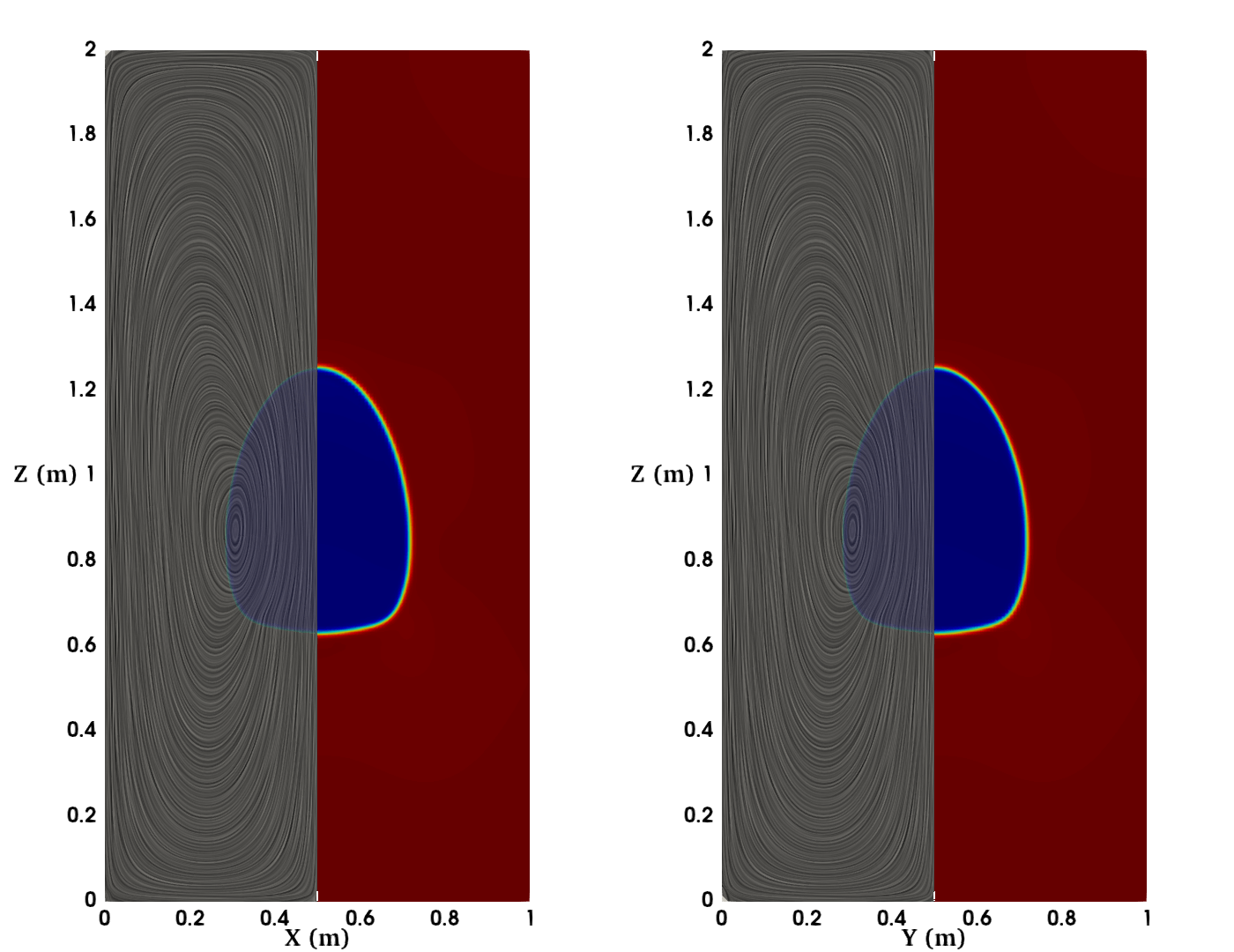}
    }
    \subfloat[$t^*=3$]{
        \includegraphics[width=0.3\linewidth]{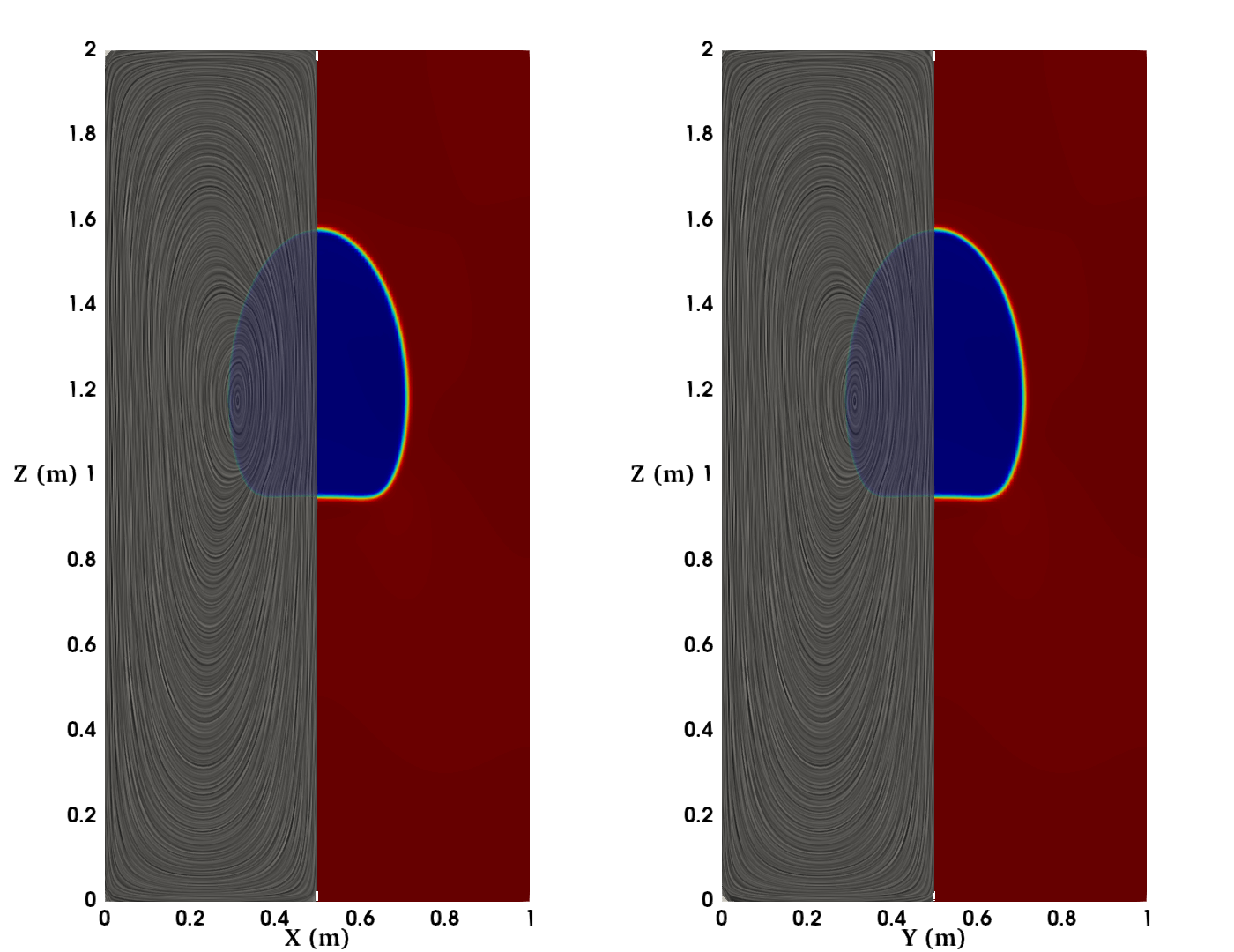}
    }
    \\
    \subfloat[$t^*=1$]{
        \includegraphics[width=0.3\linewidth]{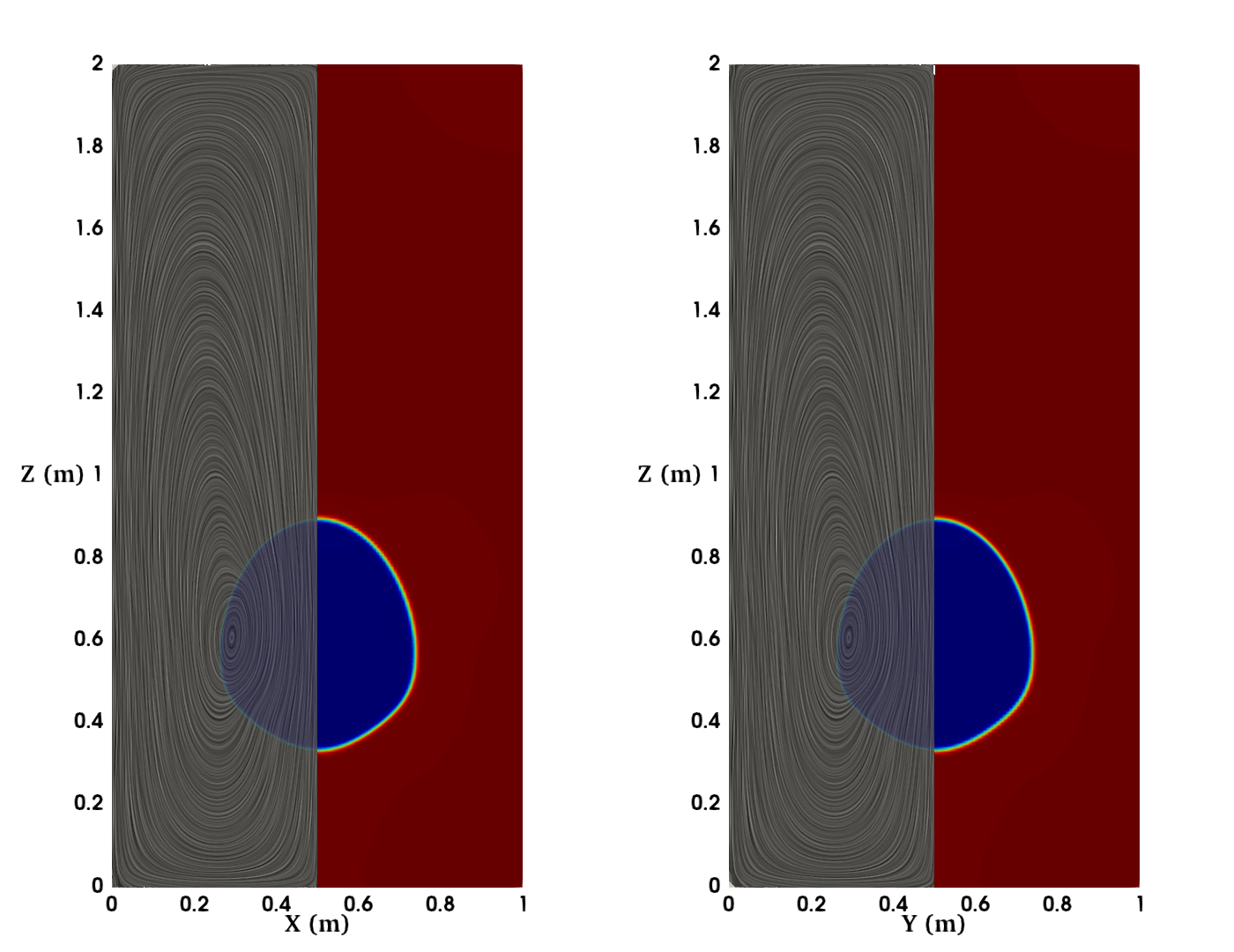}
    }
    \subfloat[$t^*=2$]{
        \includegraphics[width=0.3\linewidth]{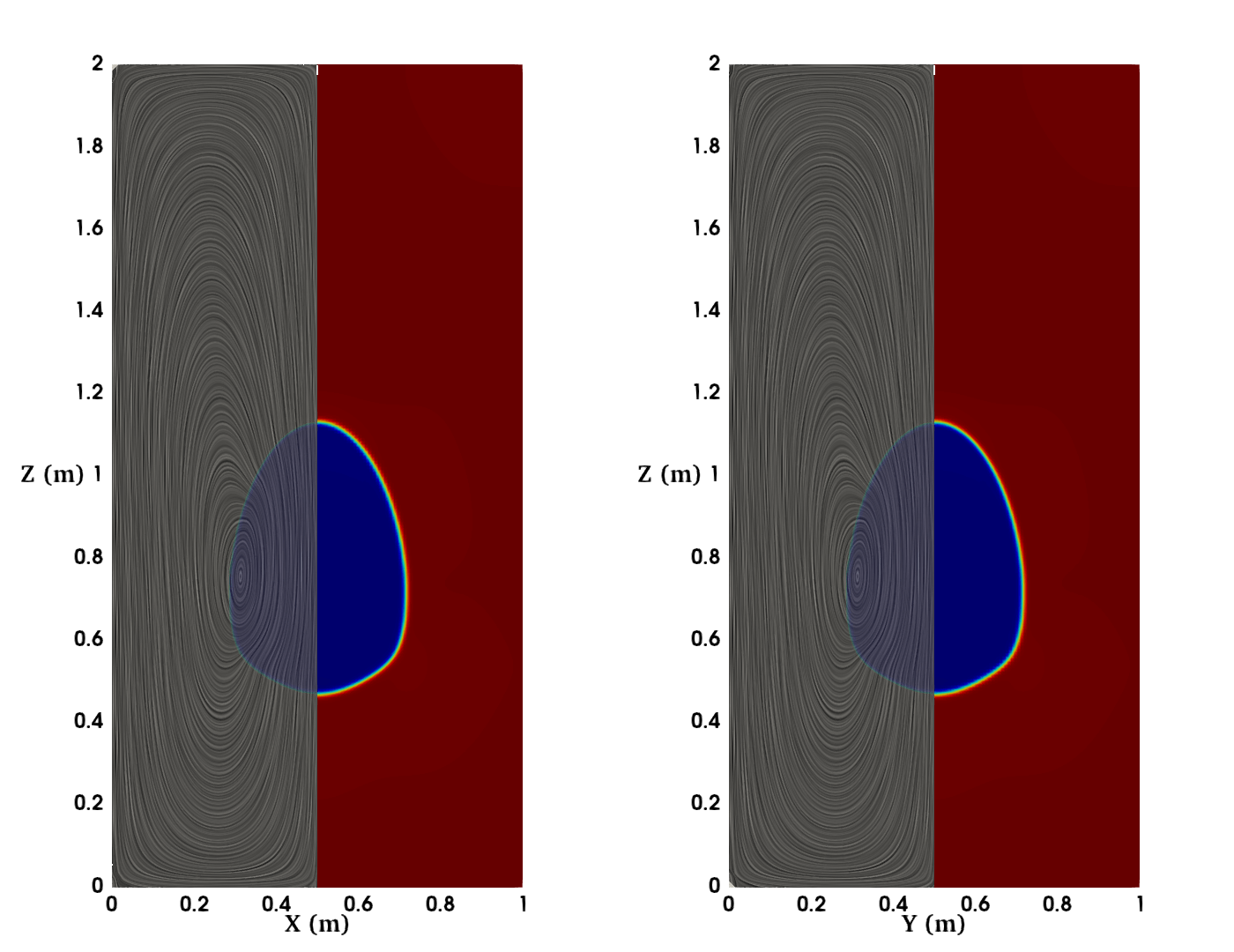}
    }
    \subfloat[$t^*=3$]{
        \includegraphics[width=0.3\linewidth]{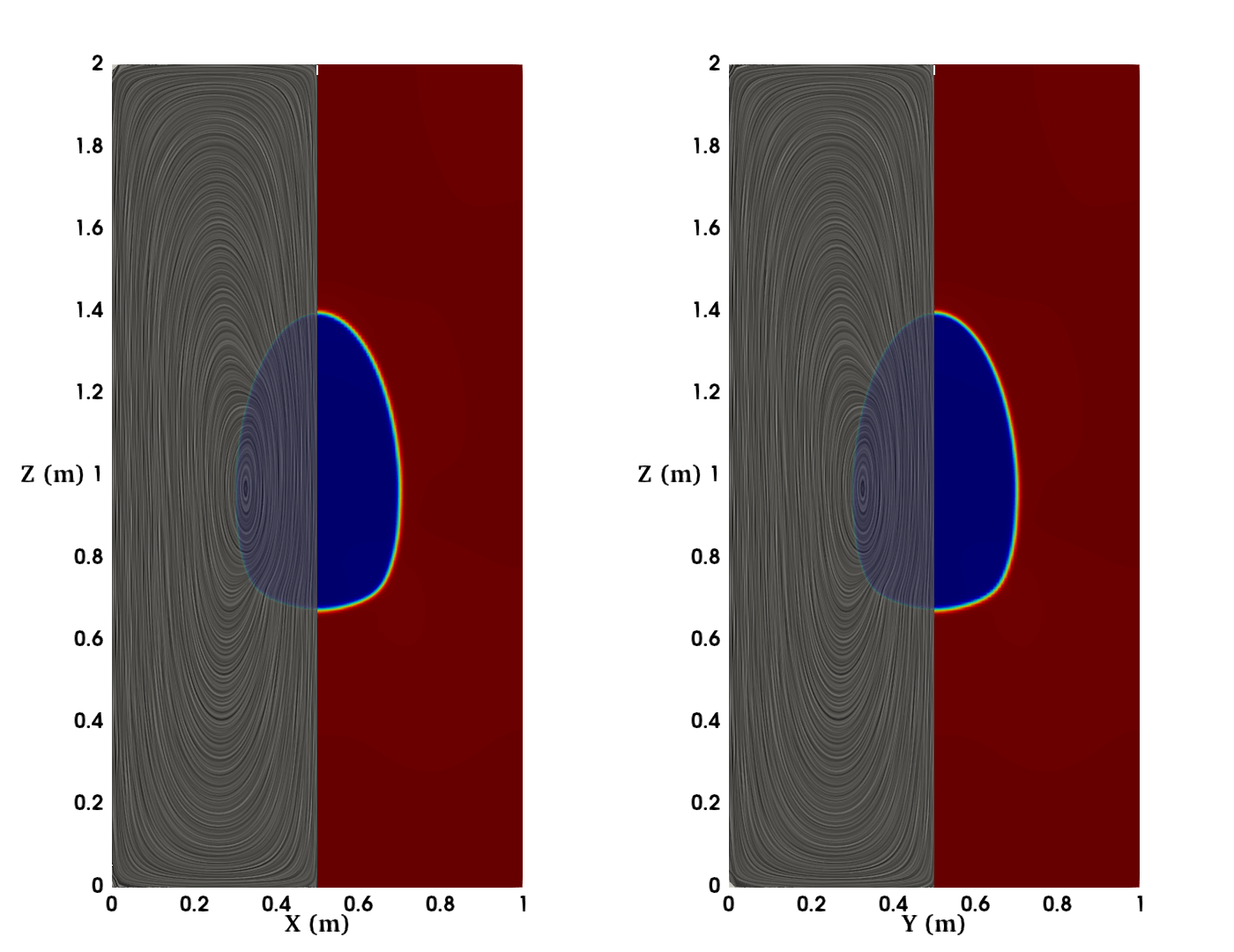}
    }
    \caption{Cross-sections of $\phi_h^n$ under vertical magnetic fields, with streamline distributions of $\bm{u}_h^n$ depicted on the left-hand side. Top row: $B_r=3$; middle row: $B_r=5$; bottom row: $B_r=7$.}
    \label{bubble_ver}
\end{figure}

\begin{figure}[H]
    \centering
    \subfloat[Horizontal magnetic fields]{
        \includegraphics[width=0.6\linewidth]{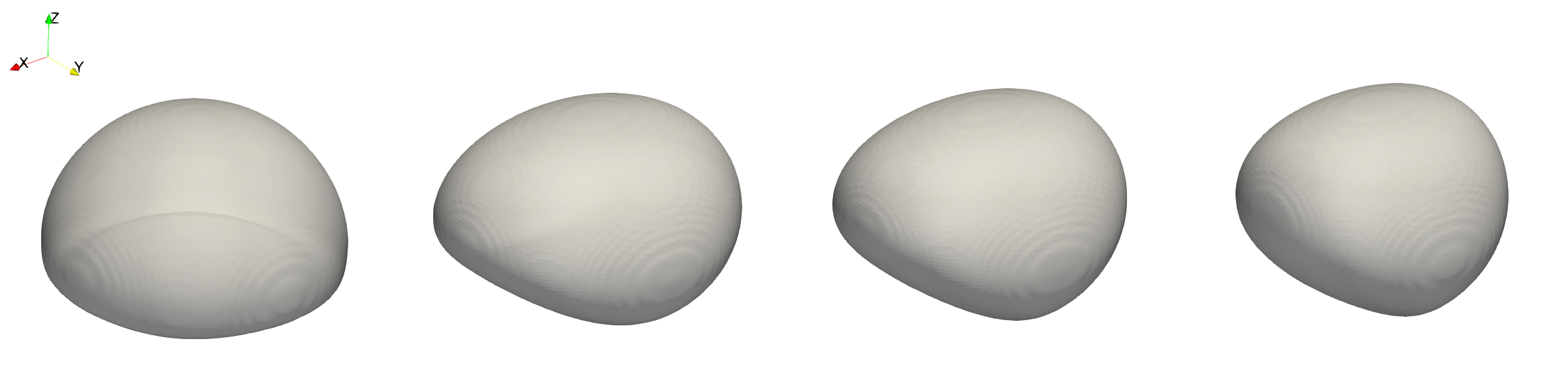}
    }
    \vspace{0.1in}\\
    \subfloat[Vertical magnetic fields]{
        \includegraphics[width=0.6\linewidth]{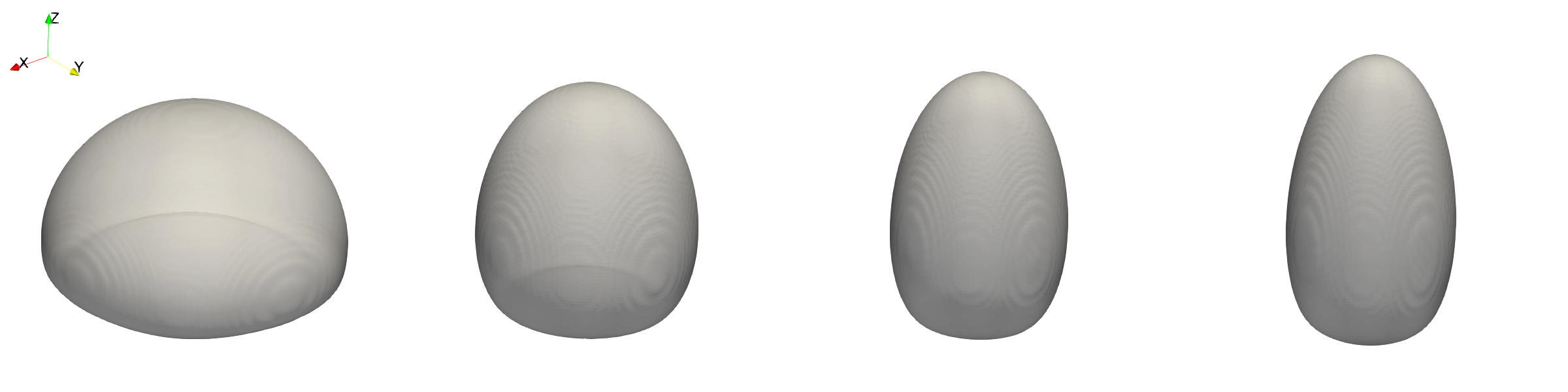}
    }
    \caption{Final bubble shape. Left to right: $B_r=0$, $3$, $5$, and $7$.}
    \label{bubble_com}
\end{figure}

The cross-sections of the order parameter $\phi$ at different dimensional time instants $t^*$ are presented in Figure \ref{bubble_ref} (in the absence of magnetic fields), Figure \ref{bubble_hor} (under horizontal magnetic fields), and Figure \ref{bubble_ver} (under vertical magnetic fields), and a comparison of final bubble shape can be found in Figure \ref{bubble_com}. The corresponding velocity streamlines, visualized via the ``Surface LIC" configuration in \texttt{ParaView}, are also depicted in these figures.

\begin{figure}[H]
	\centering
	\subfloat[Horizontal magnetic fields]{
		\includegraphics[width=0.8\linewidth]{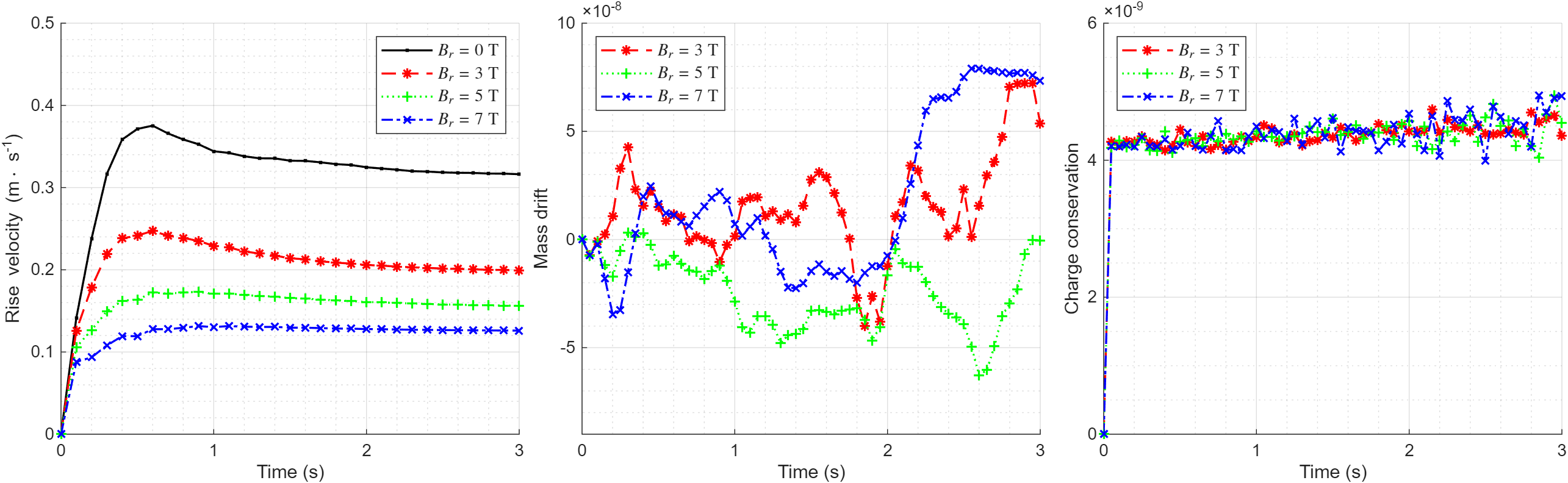}
	}
	\vspace{0.1in}\\
	\subfloat[Vertical magnetic fields]{
		\includegraphics[width=0.8\linewidth]{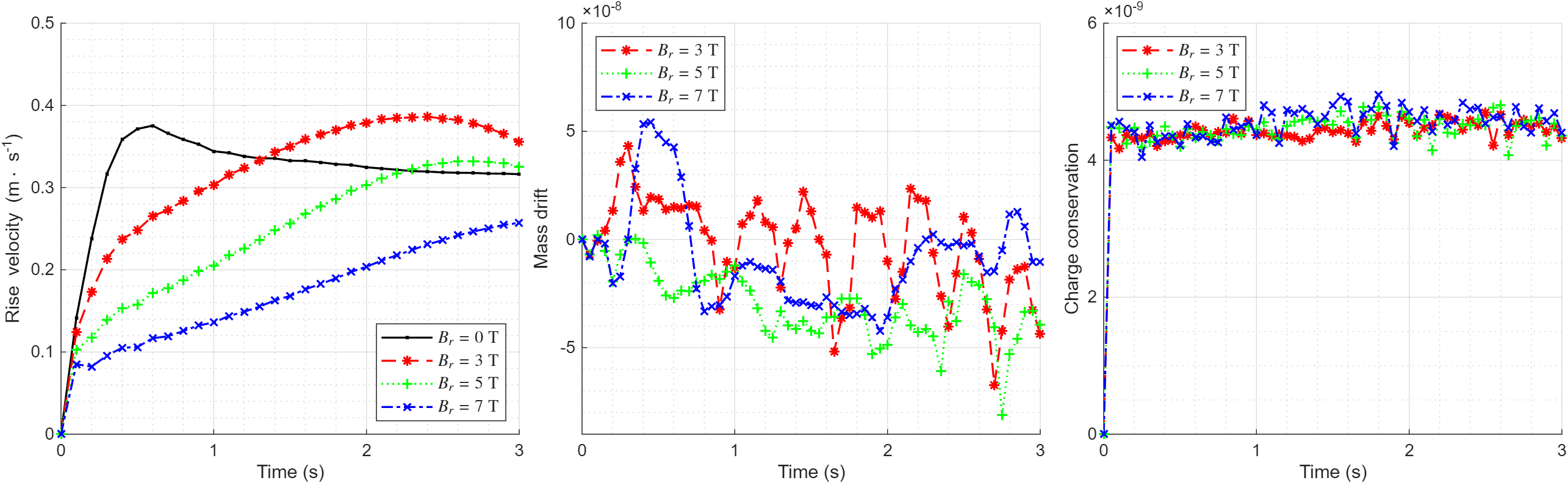}
	}
	\caption{Derived benchmark quantities.}
	\label{bubble_ben}
\end{figure}

In addition, Figure \ref{bubble_ben} shows the rise velocity, mass drift, and charge conservation, which are defined as
\begin{equation*}
    \frac{\int_{\phi_h^n<0}\bm{u}_h^n\cdot\bm{e}_3\mathrm{d}\bm{x}}{\int_{\phi_h^n<0}\mathrm{d}\bm{x}},
    \qquad
    \int_\Omega\phi_h^n\mathrm{d}\bm{x}-\int_\Omega\phi_h^0\mathrm{d}\bm{x},
    \qquad \mbox{and} \qquad
    \|\nabla\cdot\bm{J}_h^n\|_{L^2(\Omega)},
\end{equation*}
respectively. For all the reported cases, it is clearly observable that the change in the total mass and the conservation of charge is of the order $\mathcal{O}(10^{-8})$ and $\mathcal{O}(10^{-9})$. This illustrates the good mass and charge conservation property of our numerical method in practical physical simulations.

Moreover, the following conclusions can be drawn:
\begin{itemize}
    \item \textbf{velocity characteristics:} The horizontal magnetic fields exert a significant damping effect on the bubble rise velocity. In contrast, the vertical fields display a non-monotonic influence. As the magnetic intensity increases, the velocity field becomes increasingly parallel with the magnetic fields, which is clearly demonstrated by the streamline distributions.
    \item \textbf{interface dynamics:} The vertical magnetic fields preserve an axisymmetric bubble shape similar to non-magnetic scenarios. Conversely, the horizontal magnetic fields give rise to anisotropic deformations, compressing the bubble shape along the field direction. And in all cases, a stronger magnetic field leads to a progressive suppression of bubble motion.
\end{itemize}
These numerical observations are in qualitative agreement with previous results reported in \cite{2014_Zhang, 2014_Zhang_pof, 2016_Zhang}, in turn validating the capacity of our model to capture the magnetic damping effects on bubble dynamics.

\section{Concluding remarks}
\label{section_4}

In this paper, we developed a novel diffuse interface model to simulate inductionless MHD free surface flows. By means of Onsager's variational principle and the laws of thermodynamics, we derived a nonlinear system that couples the Cahn-Hilliard equation with the inductionless MHD equations. This model is thermodynamically consistent and satisfies fundamental conservation laws. Meanwhile, compared to existing diffuse interface MHD models, it is capable of handling general material properties in practical engineering processes. We carried out an asymptotic analysis via the method of formally matched asymptotic expansions. Provided a scaling law for the mobility from \cite{2013_Magaletti}, we deduced that the classical sharp interface system \eqref{sharp_interface} can be recovered in the vanishing interface thickness limit, which illustrates the robustness of our proposed diffuse interface model as an approximate approach. Regarding numerical implementation, we design an efficient decoupled, linear, and charge-conservative finite element algorithm based on a dimensionless form of the diffuse interface model. Employing this algorithm, we investigated the magnetic damping effects through a single rising bubble benchmark. The observed flow mechanisms have been observed to be consistent with the results reported in \cite{2014_Zhang, 2014_Zhang_pof, 2016_Zhang}, in turn validating the proposed diffuse interface method.

Although the developed method is capable of effectively capturing MHD phenomena in free surface flows, such as magnetic damping effects, our discussions are still relatively circumscribed. We identify several directions for future research. First, it is essential to formulate a general system that incorporates the complete MHD equations. Meanwhile, taking thermocapillary effects into account is also indispensable, in view of practical industrial applications \cite{2018_Thomas, 2014_Zhang}. Another vital aspect for future investigation lies in the exploration of other MHD phenomena within the proposed method, for instance, metal pad rolling in multilayer liquid metal systems \cite{2003_Gerbeau, 2019_Herreman, 2018_Tucs}, which further involves moving contact line dynamics and exhibits multiscale features. Finally, it is highly desirable to develop structure- and asymptotic-preserving numerical schemes that can preserve fundamental physical properties and obtain more accurate simulation results through convenient formulations.


\section*{Acknowledgments}
The authors would like to thank Prof. Wenyu Lei at UESTC for the fruitful discussions on numerical implementations using \texttt{deal.II}.  The research of MJL is supported in part by the NSFC 12271082, 62231016. The research of ZYX is supported in part by the NSFC 12471371. The research of LWX is supported in part by the NSFC 12431015, 62231016.

\small
\bibliographystyle{abbrv}
\bibliography{ref}
\end{document}